\documentclass[11pt]{article}
  \usepackage[letterpaper, left=1in, right=1in, top=1in, bottom=1in]{geometry}
  
  \usepackage[affil-it]{authblk}  

  \usepackage{palatino}  
  \usepackage{setspace}  
  \usepackage{color}  
  \usepackage[dvipsnames]{xcolor}  

	\usepackage[font=small,labelfont=bf]{caption}	
	\usepackage{subcaption}

  \usepackage{enumitem}  

\makeatletter
\newcommand*{\textoverline}[1]{$\overline{\hbox{#1}}\m@th$}
\makeatother

  \usepackage{amsmath}  
  \usepackage{amssymb}  
  \usepackage{amsthm}  
  \usepackage{bm}  
  \usepackage{bbm}  
  \allowdisplaybreaks[4]  
	\usepackage{relsize}  
	\usepackage{algorithm}  
	\usepackage{algpseudocode}  
	\usepackage{mathtools} 
  
  \DeclareMathOperator*{\minimize}{minimize}
  \DeclareMathOperator*{\maximize}{maximize}

  \newcommand{\sto}{\mathrm{subject\;to}}

	\theoremstyle{definition}

	\newtheorem{remark}{Remark}
	
  \usepackage{pifont}  

	\usepackage{longtable}  
	\usepackage{multirow}  
	\usepackage{array}  
	\usepackage{booktabs}  
    
  \usepackage{graphicx}  
  \usepackage{float}  
	\usepackage{hyperref}
	\hypersetup{
		pdfstartview = {FitV},  
		colorlinks = true,    
		linkcolor = NavyBlue,  
		citecolor = NavyBlue,  
		filecolor = NavyBlue,  
		urlcolor = NavyBlue  
	}  

\usepackage[numbers, sort&compress]{natbib}  

\usepackage{threeparttable}

\usepackage{bigints}

\usepackage{tikz}
\usetikzlibrary{decorations.pathreplacing,calc}
\newcommand{\tikzmark}[1]{\tikz[overlay,remember picture] \node (#1) {};}

\newcommand*{\AddNote}[4]{%
    \begin{tikzpicture}[overlay, remember picture]
        \draw [decoration={brace,amplitude=0.5em},decorate,ultra thick,black]
            ($(#3)!(#1.north)!($(#3)-(0,1)$)$) --  
            ($(#3)!(#2.south)!($(#3)-(0,1)$)$)
                node [align=center, text width=6cm, pos=0.5, anchor=west] {#4};
    \end{tikzpicture}
}%

\newcommand*{\rom}[1]{\expandafter\@slowromancap\romannumeral #1@}

\title{Column generation for multistage stochastic mixed-integer nonlinear programs with discrete state variables}
\author[1]{Tushar Rathi}
\author[1]{Benjamin P. Riley}
\author[2,3]{Angela Flores‑Quiroz}
\author[1]{Qi Zhang \thanks{Corresponding author (qizh@umn.edu)}}
\affil[1]{Department of Chemical Engineering and Materials Science, University of Minnesota, Minneapolis, MN 55455, USA}
\affil[2]{Department of Electrical Engineering, University of Chile, Santiago, Chile}
\affil[3]{Instituto Sistemas Complejos de Ingeniería, Santiago, Chile}
\date{}

\begin{document}

\maketitle

\begin{abstract}
\noindent  Stochastic programming provides a natural framework for modeling sequential optimization problems under uncertainty; however, the efficient solution of large-scale multistage stochastic programs remains a challenge, especially in the presence of discrete decisions and nonlinearities. In this work, we consider multistage stochastic mixed-integer nonlinear programs (MINLPs) with discrete state variables, which exhibit a decomposable structure that allows its solution using a column generation approach. Following a Dantzig-Wolfe reformulation, we apply column generation such that each pricing subproblem is an MINLP of much smaller size, making it more amenable to global MINLP solvers. We further propose a method for generating additional columns that satisfy the nonanticipativity constraints, leading to significantly improved convergence and optimal or near-optimal solutions for many large-scale instances in a reasonable computation time. The effectiveness of the tailored column generation algorithm is demonstrated via computational case studies on a multistage blending problem and a problem involving the routing of mobile generators in a power distribution network.

\bigbreak


 \end{abstract}

\section{Introduction}
Decision-making problems under uncertainty are often formulated using the stochastic programming \citep{birge2011introduction} framework, especially if the distribution of the uncertain parameters is known a priori. Consider a general multistage stochastic programming problem of the following form:
\begin{align*}
    \minimize_{x_{1},y_{1} \in \mathcal{G}_{1}(\xi_{1})} f_{1}(x_{1},y_{1},\xi_{1}) + \mathbb{E}\Big( \min_{x_{2},y_{2} \in \mathcal{G}_2(x_{1},\xi_{1})} f_{2}(x_{2},y_{2},\xi_{[2]}) + \mathbb{E} \Big( \cdots + \\
    \mathbb{E} \Big( \min_{x_{T},y_{T} \in \mathcal{G}_{T}(x_{T-1},\xi_{[T]})} f_{T}(x_{T},y_{T},\xi_{[T]}) |\xi_{[T-1]}\Big) \cdots |\xi_{[2]}\Big)|\xi_{1}\Big),
\end{align*}
where $x$ and $y$ are the state and stage variables, respectively. The state variables link decisions in different stages, whereas the stage variables are local to a stage. All uncertainties realized up to stage $t$ are represented by $\xi_{[t]}$, i.e. $\xi_{[t]} =(\xi_1, \dots, \xi_t)$. The above formulation enables determining optimal decisions in a sequential decision-making process that anticipates uncertainty in parameters that realizes over time.

Multistage stochastic programming is a well-studied framework that has seen significant theoretical advances in both modeling and solution methodologies. It is frequently used to model problems in applications with high uncertainty, such as financial planning, long-term expansion planning, power systems operations, inventory management, and supply chain engineering. A sample average approximation approach \citep{kleywegt2002sample} is commonly applied, resulting in stochastic programming models that scale with the number of discrete scenarios considered, which can be particularly large in the multistage setting. Various solution methods have been developed over the years to improve the tractability of such stochastic programs. However, most existing methods are targeted at solving two-stage or multistage linear stochastic programs. For a comprehensive review of decomposition techniques in stochastic optimization, we refer to \citet{escudero2017scenario}. In this work, our focus is on multistage stochastic programs with discrete state variables and generally nonlinear problems. In the following, we review works focused on multistage stochastic programs with discrete decisions and\slash or those involving nonlinearities.

In the case of problems involving discrete variables, significant algorithmic advancements have been made for two-stage stochastic mixed-integer programs. However, improving the tractability of multistage stochastic programs with discrete decisions using decomposition-based algorithms remains an active research area. \citet{zou2019stochastic} extended the idea of stochastic dual dynamic programming (SDDP) \citep{pereira1991multi, shapiro2011analysis} by proposing stochastic dual dynamic integer programming (SDDiP), which can handle binary state variables in multistage stochastic programs. \citet{lara2020electric} extended it further to problems with mixed-integer state variables, although without guaranteed finite convergence due to potential duality gaps. Recently, \citet{ahmed2022stochastic} proposed stochastic Lipschitz dynamic programming for solving multistage stochastic mixed-integer programs. Note that \citet{zou2019stochastic}, \citet{lara2020electric}, and \citet{ahmed2022stochastic} assume stagewise independence of uncertainty, whereas the framework we propose in this paper does not require this assumption. The progressive hedging algorithm has also been proposed as a heuristic for solving stochastic mixed-integer programming problems \citep{rockafellar1991scenarios,lokketangen1996progressive}. However, even for the two-stage case, convergence remains an issue \citep{watson2011progressive}. \citet{barnett2017bbph} and \citet{atakan2018progressive} integrated progressive hedging within the branch-and-bound algorithm. While these approaches led to reduced optimality gaps, they also resulted in increased solution times.  The column generation algorithm has emerged as a promising technique for enhancing tractability \citep{singh2009dantzig, sen2006stochastic, flores2021distributed, rathi2022capacity} of multistage stochastic mixed-integer linear programming (MILP) problems. \citet{flores2021distributed} showed its superiority over progressive hedging and nested Benders decomposition for power system planning under uncertainty, proposing additional convergence enhancements through the concept of column sharing. Although the aforementioned papers \citep{singh2009dantzig, sen2006stochastic, flores2021distributed, rathi2022capacity} demonstrate the use of column generation for solving stochastic MILPs, its application and potential for solving stochastic mixed-integer nonlinear programs (MINLPs) have been largely overlooked in the literature, which we address in this paper.

The computational complexity further escalates when, in addition to discrete decisions, nonlinearities are present in the model. From the perspective of dealing with nonlinearity, most efforts have focused on two-stage stochastic problems. \citet{li2011nonconvex} proposed a nonconvex generalized Benders decomposition algorithm for solving two-stage stochastic MINLPs with binary first-stage and continuous recourse variables. \citet{cao2019scalable} proposed a reduced-space branch-and-bound scheme that solves lower- and upper-bounding problems for two-stage stochastic nonlinear programs by relaxing the nonanticipativity constraints and solving scenario problems (resulting from fixing first-stage decisions), respectively. \citet{li2018improved} dealt with two-stage convex mixed-binary variables in both stages via an improved L-shaped method with Lagrangean and strengthened Benders cuts in the Benders master problem, followed by developing a generalized Benders decomposition-based branch-and-bound algorithm with finite $\epsilon-$convergence for the same class of problems \citep{li2019finite}. They further proposed a generalized Benders decomposition-based branch-and-cut approach for small-sized nonconvex two-stage stochastic problems with mixed-binary variables in both stages \citep{li2019generalized}. \citet{allman2021branch} proposed a branch-and-price algorithm to solve a class of nonconvex MINLPs with integer linking variables. There have been very limited efforts in developing efficient solution methods for multistage stochastic MINLPs. Recently, \citet{zhang2022stochastic} proposed a generalization of the SDDP approach to multistage stochastic MINLPs. Their method incorporates a cut generation scheme and derives a surrogate representation of the original model through a regularization approach. \citet{fullner2022non} developed a nonconvex nested Benders decomposition algorithm for multistage MINLP problems; in theory, this method can solve multistage stochastic MINLPs but the paper has, however, only demonstrated its application to deterministic problems.

In summary, the majority of the research on algorithmic advances for multistage stochastic programs has focused on MILPs, while mostly two-stage stochastic programs have been considered in regard to MINLPs. Hence, given that multistage stochastic MINLPs have received far less attention despite their ability to model various real-world problems, in this work, we propose a decomposition strategy based on column generation targeting multistage stochastic MINLPs featuring discrete state variables. In the following, we summarize our main contributions.

\begin{enumerate}
\item We develop a Dantzig-Wolfe reformulation of the general multistage stochastic MINLP via a discretization approach, which enables the decomposition of the problem using column generation. This method confines the nonlinearity to smaller-sized subproblems, making them amenable to off-the-shelf MINLP solvers.

\item While each subproblem can be solved independently, enabling parallelization and leading to reduced solution times, we additionally adopt the column sharing strategy from \citet{flores2021distributed} that exploits the scenario tree structure to facilitate information (column) sharing among subproblems (nodes in the scenario tree). This improves convergence by providing more information to the master problem of the column generation algorithm.

\item We showcase the efficacy of our decomposition scheme via two computational case studies: one on multistage blending and the other on the routing of mobile generators in a power distribution network. In particular, we conduct a comparative analysis of the proposed column generation approach, both with and without column sharing, against solutions obtained from directly solving the fullspace model. We also provide insights into how incorporating column sharing into the column generation algorithm can lead to improved performance.

\end{enumerate}

The remainder of this paper is organized as follows. In Section \ref{sec:reformulation}, we discuss the Dantzig-Wolfe reformulation of the multistage stochastic MINLP formulation with discrete state variables, followed by the column generation decomposition scheme in Section \ref{sec:CG}. Section \ref{sec:CS} describes the column sharing algorithm and how it may lead to the generation of better columns in each iteration. Sections \ref{sec:blending} and \ref{sec:OPF} focus on the case studies on multistage blending and the routing of mobile generators in a power distribution network, respectively. Finally, we provide concluding remarks in Section \ref{sec:conclusions}.

\section{Reformulation via discretization}
\label{sec:reformulation}
\begin{figure}[H]
\centering    \includegraphics[width=0.7\linewidth]{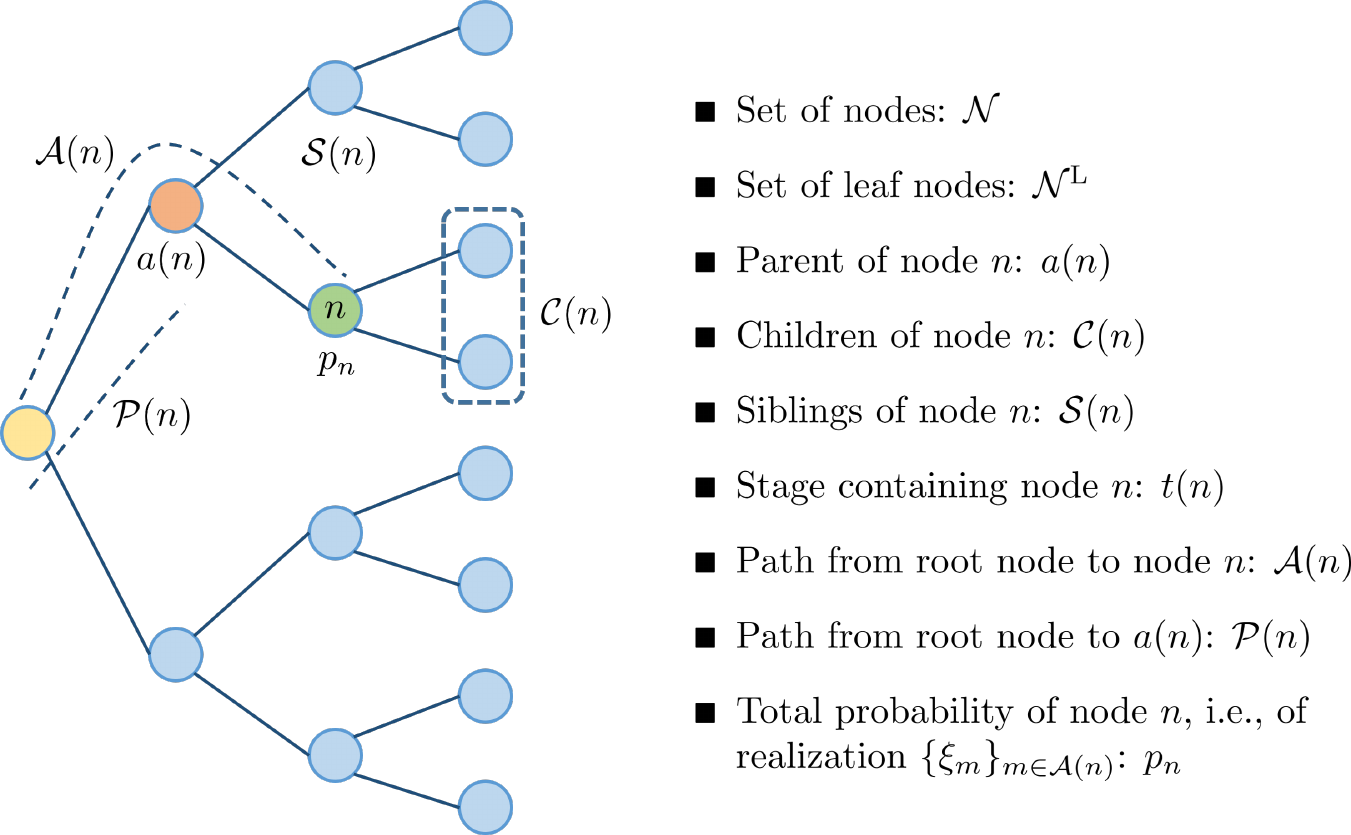}
    \caption{Schematic of a typical scenario tree for a multistage stochastic programming problem.}
    \label{fig:scenario_tree}
\end{figure}

Consider the following extensive form of a general multistage stochastic programming formulation with discrete state variables:
\begin{align*}
\text{(SP)} \qquad \minimize_{x,y} \quad & \sum_{n \in \mathcal{N}} p_{n} f_{t(n)}(x_{n}, y_{n}, \xi_{m|m \in \mathcal{A}(n)})\\
    \sto \quad & g_{t(n)}(x_{a(n)}, x_{n}, y_{n}, \xi_{m|m \in \mathcal{A}(n)}) \leq 0 \quad \forall \, n \in \mathcal{N} \\
    & x_{n} \in \mathbb{Z}^{d_{t(n)}}, y_{n} \in \mathbb{R}^{q_{t(n)}} \times \mathbb{Z}^{r_{t(n)}} \quad \forall \, n \in \mathcal{N},
\end{align*}
where the set of nodes in the scenario tree (see Figure \ref{fig:scenario_tree}) is denoted by $\mathcal{N}$, and the probability corresponding to a node $n \in \mathcal{N}$ is $p_{n}$. The parent node of a node $n$ is denoted by $a(n)$, and the set of ancestor nodes of node $n$, i.e. the nodes on the path connecting the root node to node $n$, is denoted by $\mathcal{A}(n)$. The cost function for stage $t$, $f_{t(n)}$, applies to each node $n$ in that particular stage and is a function of the state variables $x_{n}$, stage variables $y_{n}$, and the uncertainty realized up to that point in time $\xi_{m|m \in \mathcal{A}(n)}$. Similarly, $g_{t(n)}$ represents stage-specific constraints that may be functions of the state variables from the previous stage, the current state and stage variables, and the uncertainties realized up to that stage.

In this work, we consider the following specific form of (SP), which constitutes a common class of multistage stochastic MINLPs:
\begin{align}
\text{(MSSP)} \qquad z_{\mathrm{MSSP}} = \min_{x,y,u} \;\; \quad \quad &  c^\top u + \sum_{n \in \mathcal{N}} p_{n}f_{t(n)}(x_{n}, y_{n})\\
    \sto \quad 
    & Au + \sum_{n \in \mathcal{N}}D_nx_{n} \geq b \label{complicating_constraint} \\
    & g_{t(n)}(x_{n}, y_{n}) \leq 0 \quad \forall \, n \in \mathcal{N} \label{stage_specific_constraint} \\
    & x^{\mathrm{min}}_{t(n)} \leq x_{n} \leq x^{\mathrm{max}}_{t(n)} \quad \forall \, n \in \mathcal{N} \label{state_var_bounds} \\
    & u \in \mathbb{R}^{m} \times \mathbb{Z}^{\bar{m}} \label{bound_u} \\
    & x_{n} \in \mathbb{Z}^{d_{t(n)}}, y_{n} \in \mathbb{R}^{q_{t(n)}} \times \mathbb{Z}^{r_{t(n)}}\quad \forall \, n \in \mathcal{N}. 
\end{align}
Here, for brevity, we omit $\xi_{m|m \in \mathcal{A}(n)}$, but it must be understood that the equations defining the constraints and objective terms at node $n$ are functions of the uncertainty realized up to that point in time. State variables, $x_{n}$, are restricted to integer values and are bounded according to constraints \eqref{state_var_bounds}, whereas stage variables, $y_{n}$, can take integer as well as continuous values. Additionally, we may have variables $u$ that appear in constraints with only state variables. Constraints \eqref{complicating_constraint} link the decisions across different time periods via the state variables whereas constraints \eqref{stage_specific_constraint} are stage-specific. Functions $f_{t(n)}(x_{n}, y_{n})$ and $g_{t(n)}(x_{n}, y_{n})$ can be generally nonlinear and nonconvex. It is worth noting that, if not already the case, a stochastic program can often be transformed into the above form using auxiliary variables. An example of such a reformulation using auxiliary variables is presented later in our case study on the multistage blending problem.

Now, any multistage stochastic programming model with this structure can be decomposed into $|\mathcal{N}|$ subproblems, provided we remove the complicating constraints \eqref{complicating_constraint}. For each subproblem associated with node $n \in \mathcal{N}$, the feasible set for the state variables $x_{n}$ can be defined as follows:
\begin{align}
    \mathcal{X}_{n} := \{x_{n} \in \mathbb{Z}^{d_{t(n)}}: \exists \,\, y_{n} \in \mathbb{R}^{q_{t(n)}} \times \mathbb{Z}^{r_{t(n)}} \,\, \text{such that} \,\, g_{t(n)}(x_{n}, y_{n}) \leq 0, \,\, x^{\mathrm{min}}_{t(n)} \leq x_{n} \leq x^{\mathrm{max}}_{t(n)} \}.
\end{align}

Since $x_{n}$ are integer and bounded, $\mathcal{X}_{n}$ is a finite set. Accordingly, we can rewrite $\mathcal{X}_{n} = \{x_{n1}^{*}, \ldots, x_{n,K_{n}}^{*} \}$, where $|\mathcal{X}_{n}| = K_{n}$ and $x_{nk}^{*}$ is the $k^{th}$ feasible $x_{n}$. However, in the final solution, $x_{n}$ takes exactly one of these feasible values, which can be enforced by the following constraints:
\begin{align}
& x_{n} = \sum_{k \in \mathcal{K}_{n}} \rho_{nk}x_{nk}^{*} \\
& \sum_{k \in \mathcal{K}_{n}} \rho_{nk} = 1 \\ 
& \rho_{nk} \in \{0,1\} \quad \forall \, k \in \mathcal{K}_{n},
\end{align}
where $\mathcal{K}_{n} = \{1,\ldots,K_{n}\}$. The cost associated with the column $x_{nk}^{*} \in \mathcal{X}_{n}$ can be obtained by solving the following optimization problem:
\begin{align} 
    f_{nk}^{*} = \min_{y_{n} \in \mathbb{R}^{q_{t(n)}} \times \mathbb{Z}^{r_{t(n)}} } \{ p_{n}f_{t(n)}(x_{nk}^{*},y_{n}): g_{t(n)}(x_{nk}^{*}, y_{n}) \leq 0 \}.
\end{align}

Now that we know how to select exactly one column from each subproblem's finite feasible space, as well as the associated cost, we use the discretization approach \citep{vanderbeck2005implementing, singh2009dantzig} to reformulate (MSSP) (also known as Dantzig-Wolfe reformulation) as follows:
\begin{align} \hspace{-5em}
\text{(M-MSSP)} \qquad z_{\mathrm{MSSP}}^{\mathrm{M}} = \min_{u, \rho} \;\; \quad \quad & c^\top u + \sum_{n \in \mathcal{N}}\sum_{k \in \mathcal{K}_{n}} \rho_{nk}f_{nk}^{*} \\
    \sto \quad
    & Au + \sum_{n \in \mathcal{N}}D_{n}\sum_{k \in \mathcal{K}_{n}}\rho_{nk}x_{nk}^{*} \geq b \label{MMSSP_linking} \\
    & \sum_{k \in \mathcal{K}_{n}} \rho_{nk} = 1 \quad \forall \, n \in \mathcal{N} \label{MMSSP_convexity}\\
    & u \in \mathbb{R}^{m} \times \mathbb{Z}^{\bar{m}} \quad \forall \, n \in \mathcal{N} \\
    & \rho_{nk} \in \{0,1\} \quad \forall \, n \in \mathcal{N}, k \in \mathcal{K}_{n}.
\end{align}
The two formulations, (MSSP) and (M-MSSP), are equivalent \citep{vanderbeck2005implementing} in that their optimal objective values are the same, i.e. $z_{\mathrm{MSSP}} = z_{\mathrm{MSSP}}^{\mathrm{M}}$, and there is a direct mapping between the solutions from the two problems.

\section{Column generation}
\label{sec:CG}
The goal of column generation \citep{barnhart1998branch, lubbecke2005selected} is to solve the linear programming (LP) relaxation of (M-MSSP), which we denote by (\textoverline{M-MSSP}). In general, solving (\textoverline{M-MSSP}) can be computationally difficult due to $K_{n}$ being an exponentially large number. For that reason, to start with, we populate (\textoverline{M-MSSP}) with a few feasible columns in set $\overline{\mathcal{K}}_{n}$, where $|\overline{\mathcal{K}}_{n}| \ll K_{n}$, and henceforth generate additional columns iteratively that can potentially improve the objective value using the following pricing problem:
\begin{align}
\text{(PP)} \qquad \psi \; =  \quad \min_{x,y} \;\; \quad \quad & \sum_{n \in \mathcal{N}} \big[p_{n} f_{t(n)}(x_{n}, y_{n}) - \gamma^{\top}D_{n}x_{n}-\mu_{n} \big] \\
    \sto \quad 
    & g_{t(n)}(x_{n}, y_{n}) \leq 0 \quad \forall \, n \in \mathcal{N} \\
    & x^{\mathrm{min}}_{t(n)} \leq x_{n} \leq x^{\mathrm{max}}_{t(n)} \quad \forall \, n \in \mathcal{N} \\
    & x_{n} \in \mathbb{Z}^{d_{t(n)}}, y_{n} \in \mathbb{R}^{q_{t(n)}} \times \mathbb{Z}^{r_{t(n)}} \quad \forall \, n \in \mathcal{N}, 
\end{align}
where $\gamma$ and $\mu$ are the dual prices corresponding to constraints \eqref{MMSSP_linking} and \eqref{MMSSP_convexity}, respectively. The idea behind solving (PP) is to find the column with the most negative reduced cost and add it to (\textoverline{M-MSSP}). In the case of a minimization problem, columns with negative reduced costs are potential candidates for improving the incumbent objective value. In general, (PP) can be nonconvex and hence difficult to solve; however, it can be further decomposed into $|\mathcal{N}|$ independent pricing problems, one for each node $n \in \mathcal{N}$. Apart from the obvious advantage of reduced solution time via parallelization, smaller-sized pricing problems are more computationally tractable when solved using commercial global solvers such as Gurobi and BARON. The pricing problem for node $n$ is as follows:
\begin{align} \hspace{-4em}
\text{(PP\textsubscript{$n$})} \qquad \psi_{n} \; = \quad \min_{x_{n},y_{n}} \;\; \quad \quad & p_{n} f_{t(n)}(x_{n}, y_{n}) - \gamma^{\top}D_{n}x_{n}-\mu_{n}\\
    \sto \quad 
    & g_{t(n)}(x_{n}, y_{n}) \leq 0 \\
    & x^{\mathrm{min}}_{t(n)} \leq x_{n} \leq x^{\mathrm{max}}_{t(n)} \\
    & x_{n} \in \mathbb{Z}^{d_{t(n)}}, y_{n} \in \mathbb{R}^{q_{t(n)}} \times \mathbb{Z}^{r_{t(n)}}.
\end{align}

The detailed steps of solving (\textoverline{M-MSSP}) via column generation are shown in Algorithm \ref{alg:CG}. Since we start with a small subset of feasible columns, as indicated in line \ref{algl:solve_relaxed_RMP} in Algorithm \ref{alg:CG}, we essentially solve a restricted version of (\textoverline{M-MSSP}), which we denote (\textoverline{RM-MSSP}). The resulting objective value, $z_{\mathrm{MSSP}}^{\overline{\mathrm{RM}}}$, from solving (\textoverline{RM-MSSP}) provides an upper bound (UB\textsuperscript{\textoverline{MP}}) to (\textoverline{M-MSSP}), i.e. $z_{\mathrm{MSSP}}^{\overline{\mathrm{M}}} \leq z_{\mathrm{MSSP}}^{\overline{\mathrm{RM}}}$. Post solving (\textoverline{RM-MSSP}) with the initial set of feasible columns $\overline{\mathcal{K}}_{n}$ from all $n \in \mathcal{N}$, as indicated in lines \ref{algl: start_solve_sp}-\ref{algl:stop_solve_sp} in Algorithm \ref{alg:CG}, in every iteration, we generate additional columns via the pricing problems (PP\textsubscript{$n$}) for all $n \in \mathcal{N}$, which also provide a lower bound (LB\textsuperscript{\textoverline{MP}}) to (\textoverline{M-MSSP}). More precisely, arguments from duality theory \citep{wolsey1998integer} can be used to show that $z_{\mathrm{MSSP}}^{\overline{\mathrm{RM}}} + \sum_{n \in \mathcal{N}}\psi_{n} \leq z_{\mathrm{MSSP}}^{\overline{\mathrm{M}}}$. The algorithm converges when the relative optimality gap drops below the desired tolerance $\epsilon$. Given that there are a finite number of possible solutions for the discrete state variables, which translate into a finite number of columns that can be potentially priced out, the algorithm is guaranteed to converge in a finite number of iterations \citep{desrosiers2005primer}.

\begin{algorithm}[htbp]
\caption{Column generation for (\textoverline{M-MSSP})}
\label{alg:CG}
\begin{algorithmic}[1]
\State $\textrm{LB\textsuperscript{\textoverline{MP}}} \gets - \infty$, $\textrm{UB\textsuperscript{\textoverline{MP}}} \gets + \infty$, 
$k_{n} \gets |\overline{\mathcal{K}}_{n}| \,\, \forall \, n \in \mathcal{N}$ \tikzmark{right} \Comment{Initialization} \label{algl:cg_initialize} 

\While{$\frac{|\textrm{UB\textsuperscript{\textoverline{MP}}} - \textrm{LB\textsuperscript{\textoverline{MP}}}|}{|\textrm{UB\textsuperscript{\textoverline{MP}}}|} > \epsilon$} 

\State Solve (\textoverline{RM-MSSP}), get $z_{\mathrm{MSSP}}^{\overline{\mathrm{RM}}}$, $u^{*}$, $\rho^{*}$, and dual prices $\gamma$ and $\mu$ \label{algl:solve_relaxed_RMP}
\State $\textrm{UB\textsuperscript{\textoverline{MP}}} \gets z_{\mathrm{MSSP}}^{\overline{\mathrm{RM}}}$ \label{algl:update_ub}
    \For{$n \in \mathcal{N}$} \tikzmark{top} \label{algl: start_solve_sp}
        \State Solve (PP\textsubscript{$n$}), get $\psi_{n}$ and $x_{n}^{*}$
        \If{$\psi_{n} < 0$}
            \State $k_{n} \gets k_{n}+1$
            \State $x_{nk_{n}}^{*} \gets x_{n}^{*}$
            \State $f_{nk_{n}}^{*} \gets \psi_{n} + \gamma^\top D_{n}x_{n}^{*}+\mu_{n}$ \label{lst:line:col_cost}
            \State $\overline{\mathcal{K}}_{n} \gets \overline{\mathcal{K}}_{n} \cup \{k_{n}\}$
        \EndIf
    \EndFor \tikzmark{bottom} \label{algl:stop_solve_sp}
\State $\textrm{LB\textsuperscript{\textoverline{MP}}} \gets \textrm{max}\{\textrm{LB\textsuperscript{\textoverline{MP}}}, z_{\mathrm{MSSP}}^{\overline{\mathrm{RM}}} + \sum_{n \in \mathcal{N}}\psi_{n} \}$ \label{algl:update_lb}
\EndWhile
\If{non-integer solution} \Comment{see Remark 2}
    \State Solve (RM-MSSP),  get $\textrm{UB\textsuperscript{MP}}$, $u^{*}$, $\rho^{*}$ \label{algl:noninteger_soln}
\EndIf
\State \textbf{return} $u^{*},  \rho^{*}$
\end{algorithmic}
\AddNote{top}{bottom}{right}{Independent pricing problems can be solved in parallel}
\end{algorithm}

\begin{remark}
Although the pricing problems are generally nonconvex MINLPs, which can be computationally difficult to solve, they do not have to be solved to optimality in each iteration \citep{allman2021branch}. A suboptimal column with a negative reduced cost is an eligible candidate to be added to (\textoverline{RM-MSSP}); however, the LB\textsuperscript{\textoverline{MP}} update needs slight modification if pricing problems are not solved to optimality. In particular, $z_{\mathrm{MSSP}}^{\overline{\mathrm{RM}}} + \sum_{n \in \mathcal{N}} \psi_{n}^{\mathrm{LB}} \leq z_{\mathrm{MSSP}}^{\overline{\mathrm{RM}}} + \sum_{n \in \mathcal{N}} \psi_{n} \leq z_{\mathrm{MSSP}}^{\overline{\mathrm{M}}}$, where $\psi_{n}^{\mathrm{LB}}$ denotes the lower bound obtained from solving (PP\textsubscript{$n$}). Additionally, the column cost $f_{nk}^{*}$ (line~\ref{lst:line:col_cost} in Algorithm \ref{alg:CG}) then equals $\psi_{n}^{\mathrm{UB}} +\gamma^\top D_{n}x_{n}^{*}+\mu_{n}$, where $\psi_{n}^{\mathrm{UB}}$ denotes the upper bound from solving (PP\textsubscript{$n$}). Although not solving pricing problems to optimality speeds up column generation, we must solve pricing problems to optimality in the last iteration to prove convergence if the final optimality gap tolerance $\epsilon \to 0$.
\end{remark}

\begin{remark}
If a non-integer solution is obtained at the end of column generation, a mixed-integer programming (MIP) version of (\textoverline{RM-MSSP}), which we denote by (RM-MSSP), can be solved to recover an integer-feasible solution (see line \ref{algl:noninteger_soln} in Algorithm \ref{alg:CG}). The only difference between (RM-MSSP) and the (\textoverline{RM-MSSP}) in the final iteration is that (RM-MSSP) also contains the integrality constraints on $\rho$, and, if applicable, on $u$. The resulting solution provides an upper bound ($\textrm{UB\textsuperscript{MP}}$) for (M-MSSP) and may already meet the desired optimality gap or be close to the optimal solution. If it doesn't meet the desired optimality gap, branch-and-price may be used \citep{vanderbeck2000dantzig}.
\end{remark}

\section{Column sharing}
\label{sec:CS}

\begin{figure}[ht]
\centering    
\includegraphics[width=0.55\linewidth]{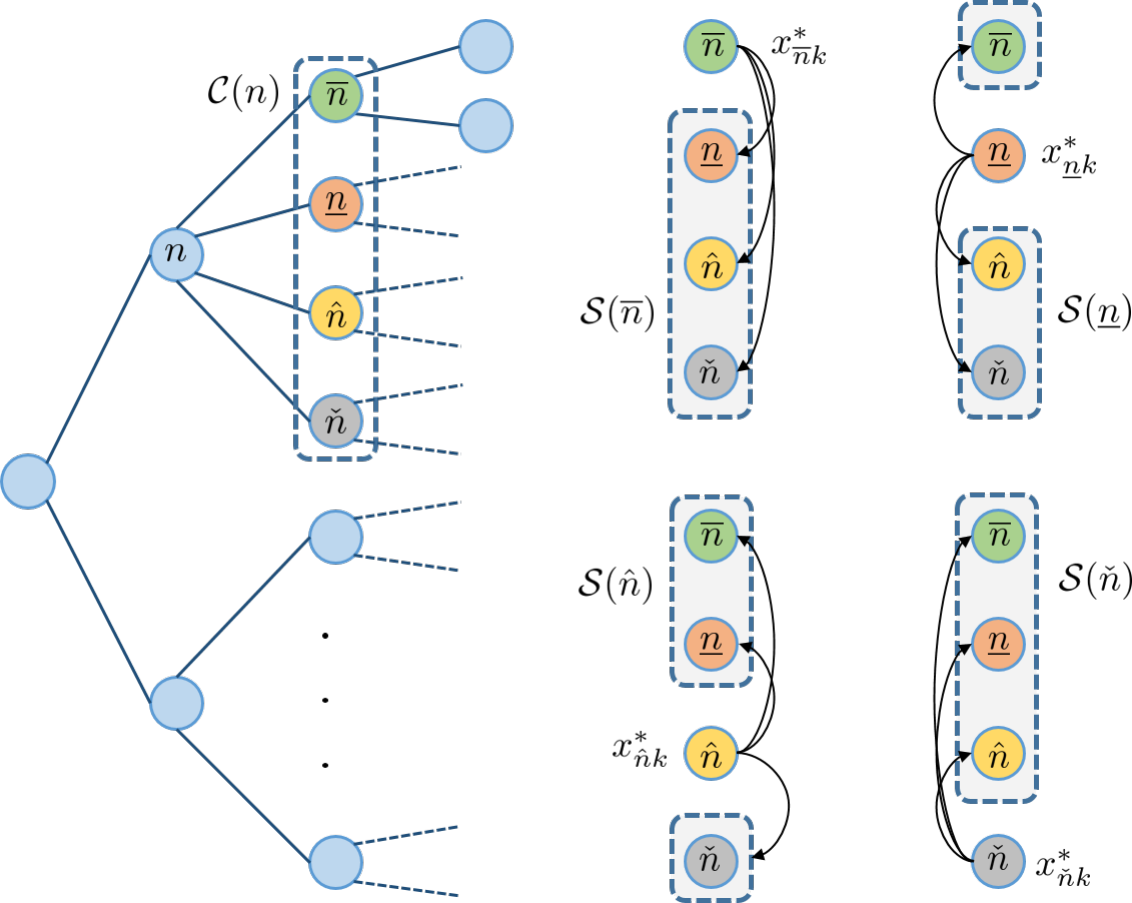}
    \caption{Illustrating column sharing amongst the first set of sibling nodes at the third stage of a scenario tree with branching structure $\mathcal{R}=\{R_1,R_2,R_3,R_4\}=\{1,2,4,2\}$, where $R_{t}$ is the number of children nodes of each node in stage $t-1$. If we price out one distinct column from each node in this scenario tree in iteration $k$, for a total of 27 distinct columns, we can obtain up to 42 additional columns by sharing them among sibling nodes.}
    \label{fig:column_sharing}
\end{figure}

Column generation is known to suffer from the \textit{heading-in} and \textit{tailing-off} effects \citep{vanderbeck2005implementing}, especially when applied to stochastic programs. One of the main reasons is that the priced out columns may be infeasible for the master problem. Because each node in the scenario tree is an independent pricing problem, the columns generated from a set of \textit{sibling nodes} (nodes with the same parent node) may not always satisfy the nonanticipativity constraints, which require having identical decisions in scenarios that are indistinguishable up to a certain point in time. Note that because in this paper we utilize the node formulation, which implicitly enforces nonanticipativity via regular model constraints, any subsequent mention of \textit{nonanticipativity constraints} is referring to the concept of nonanticipativity and not necessarily the explicit constraints that enforce it. Now, since the columns generated by sibling nodes may not necessarily satisfy nonanticipativity, this can lead to the accumulation of columns in the master problem that are often not feasible for the original problem, slowing the algorithm's convergence. To address this issue, we exploit the scenario tree structure to share columns among sibling nodes, following the approach proposed by \citet{flores2021distributed}, ensuring that a larger number of columns that satisfy the nonanticipativity constraints are generated in every iteration. The following optimization problem illustrates the sharing of a column between two sibling nodes $\overline{n}$ and $\underline{n}$:
\begin{align} \hspace{-0.5em} 
\text{(CSP\textsubscript{$\overline{n}\rightarrow\underline{n}$})} \qquad \zeta_{\overline{n} \rightarrow \underline{n} | \underline{n} \in \mathcal{S}(\overline{n})} \; = \quad \min_{x_{\underline{n}}, y_{\underline{n}}} \;\; \quad \quad & p_{\underline{n}} f_{t(\underline{n})}(x_{\underline{n}}, y_{\underline{n}})\\
    \sto \quad 
    & g_{t(\underline{n})}(x_{\underline{n}}, y_{\underline{n}}) \leq 0 \\
    & x_{\underline{n}} = x_{\overline{n}}^{*} \label{force_column} \\
    & x^{\mathrm{min}}_{t(\underline{n})} \leq x_{\underline{n}} \leq x^{\mathrm{max}}_{t(\underline{n})} \\
    & x_{\underline{n}} \in \mathbb{Z}^{d_{t(\underline{n})}}, y_{\underline{n}} \in \mathbb{R}^{q_{t(\underline{n})}} \times \mathbb{Z}^{r_{t(\underline{n})}}, 
\end{align}
where the notation $\overline{n} \rightarrow \underline{n} | \underline{n} \in \mathcal{S}(\overline{n})$ indicates that the column generated by pricing problem ($\textrm{PP}_{\overline{n}}$), $x^{*}_{\overline{n}}$, is shared with its sibling node $\underline{n}$ (constraint \
\eqref{force_column}), and $\zeta$ denotes the cost of this column with respect to node $\underline{n}$. Figure \ref{fig:column_sharing} illustrates the column sharing between a set of sibling nodes at the $k^{th}$ iteration of the column generation procedure. For a pricing problem ($\textrm{PP}_{\overline{n}}$), the column priced out can be shared with $|\mathcal{S}(\overline{n})|$ sibling nodes. Thus, assuming one column is priced out from each node, then we can generate up to $\sum_{n \in \mathcal{N}}|\mathcal{S}(n)|$ additional columns per iteration. Note that the first stage has only one node and thus undergoes no column sharing.

Here we summarize the column sharing from node $\overline{n}$ to one of its sibling node $\underline{n}$. First, the regular pricing problem for node $\overline{n}$ is solved, and the column $x^{*}_{\overline{n}}$ (with a negative reduced cost) is extracted. Now, sharing this column from node $\overline{n}$ to $\underline{n}$ implies fixing the value of that column as shown in constraint \eqref{force_column} and solving (CSP\textsubscript{$\overline{n}\rightarrow\underline{n}$}) to calculate its cost in node $\underline{n}$ (lines \ref{algl:solve_cs}-\ref{algl:cost} in Algorithm \ref{alg:column_sharing}). For large problems, columns can quickly add up, increasing the solution time of (\textoverline{RM-MSSP}). However, this could be avoided by not sharing redundant columns and discarding the infeasible ones. First, if node $\underline{n}$ has already priced out the same column with its optimal cost in the current or previous iteration, we can skip sharing it (lines \ref{algl:already_contained}-\ref{algl:skip_sharing} in Algorithm \ref{alg:column_sharing}). Second, if $x^{*}_{\overline{n}}$ is infeasible for $\underline{n}$, it can be completely discarded (lines \ref{algl:if_infeasible}-\ref{algl:exit_cs} in Algorithm \ref{alg:column_sharing}), as it will also be infeasible for the original problem. If the sharing is successful, we have created an extra column for node $\underline{n}$ in addition to its own regularly priced out column. This procedure is repeated for each pair of sibling nodes in the scenario tree. Note that sharing columns between each pair of sibling nodes is an independent problem, so we also parallelize the column sharing between all pairs of sibling nodes in each iteration.

\begin{algorithm}[H]
\caption{Column sharing from node $\overline{n}$}\label{alg:column_sharing}
\begin{algorithmic}[1]
\Procedure{Column sharing}{$\overline{n}$} \Comment{Assuming there exist a column to share}
    \State $\overline{\mathcal{N}} \gets \emptyset$ \Comment{Nodes to which column gets shared}
    \For{$\underline{n} \in \mathcal{S}(\overline{n})$} 
        \If{column from (PP\textsubscript{$\overline{n}$}) already priced out by (PP\textsubscript{$\underline{n}$}) with its optimal cost} \label{algl:already_contained}
            \State Skip sharing
        \Else \label{algl:skip_sharing}
            \State Solve (CSP\textsubscript{$\overline{n}\rightarrow\underline{n}$}) \label{algl:solve_cs}
            \If{feasible}
                \State $\overline{\mathcal{N}} \gets \overline{\mathcal{N}} \cup \{\underline{n}\} $
                \State Get $\zeta_{\overline{n} \rightarrow \underline{n}}$ \Comment{Cost of column $x_{\overline{n}}^{*}$ in node $\underline{n}$} \label{algl:cost}
            \Else \label{algl:if_infeasible}
                \State Discard column obtained from (PP\textsubscript{$\overline{n}$})
                \State \textbf{return} \Comment{If the column is infeasible, exit immediately} \label{algl:exit_cs}
            \EndIf
        \EndIf
    \EndFor 
\State \textbf{return} $\zeta_{\overline{n} \rightarrow \underline{n}} \,\, \forall \, \underline{n} \in \overline{\mathcal{N}}$
\EndProcedure
\end{algorithmic}
\end{algorithm}

\noindent The following are some key points to be aware of when applying column sharing:

\begin{enumerate}
    \item To benefit from column sharing, make sure that a shared column is required to meet the nonanticipativity constraints. If the original problem formulation does not already have a suitable structure, one can add an auxiliary variable for each node, which will refer to a common decision made at its parent node. For example, if a set of sibling nodes $\mathcal{S}(n)$ do not require $x_{n}$ to meet the nonanticipativity constraints, we can define an auxiliary variable $w_{n}:= x_{a(n)}$, where $a(n)$ represents the parent node of node $n$. Instead of $x$, we can then discretize $w$, which must satisfy nonanticipativity for all sets of sibling nodes, thereby making column sharing useful. 
    
    \item A regularly priced out column is added to (\textoverline{RM-MSSP}) only if it has a negative reduced cost; however, a column shared from node $\overline{n}$ to $\underline{n} \in \mathcal{S}(\overline{n})$, if feasible for all sibling nodes $\mathcal{S}(\overline{n})$, is added to (\textoverline{RM-MSSP}) irrespective of the sign of its reduced cost for $\underline{n}$. This is because the column has already indicated its potential to improve the objective value by having a negative reduced cost in the pricing problem for node $\overline{n}$. Sharing it with nodes in $\mathcal{S}(\overline{n})$ is primarily to increase its likelihood of becoming part of a feasible solution because it now satisfies nonanticipativity. Therefore, sharing should not be restricted even if that column has a non-negative reduced cost for any node in $\mathcal{S}(\overline{n})$.

    \item For clarity, we show sharing column from node $\overline{n}$ to $\underline{n}$ via constraint \eqref{force_column}; however, in our implementations, we simply define $x_{\underline{n}}$ as a parameter that equals $x^{*}_{\overline{n}}$, meaning that only the stage variables $y_{\underline{n}}$ remain in (CSP\textsubscript{$\overline{n}\rightarrow\underline{n}$}). This reduces the number of integer variables, making (CSP) significantly easier to solve than the regular pricing problems (PP).

    \item The goal of column sharing is to improve convergence by obtaining more feasible columns. However, in some cases, regularly priced columns may already show sufficiently good convergence; in such cases, column sharing may be used only when convergence stalls, which is typically at the beginning and end of the column generation algorithm. This may help reduce the solution time if (CSP) is computationally difficult, such as when there are a large number of integer stage variables.
\end{enumerate}

\section{Case studies}
In this section, we use two computational case studies, on multistage blending and power distribution network operation, to compare the performance of solving the fullspace multistage stochastic program, column generation (CG), and column generation with column sharing (CGCS). All models were implemented in Julia v1.9.3 \citep{Julia-2017} using the JuMP v1.20.0 \citep{Lubin2023} modeling package and solved with Gurobi 10.0.1 \citep{gurobi}. For each instance, we requested 65 cores (1 main process + 64 worker processes) to parallelize the subproblems in the CG and CGCS algorithms. The relaxed restricted master problem (\textoverline{RM-MSSP}) in both algorithms was solved using the barrier method and with presolve turned off. In most cases, column generation produced integer-feasible solutions. In a few cases where the solutions exhibited fractional values for the integer variables, we solved (RM-MSSP) to recover integer-feasible solutions, which in most instances led to solutions with gaps that either met or were very close to the optimality gap tolerance. Lastly, the number of stages in both case studies was varied by changing the number of time periods and letting uncertainty realize in each one. The number of scenarios was also varied by modifying the branching structure of the scenario tree, i.e. the number of realizations in each time period. The code for the case studies is available in our GitHub repository at \href{https://github.com/ddolab/CG-MSMINLP}{https://github.com/ddolab/CG-MSMINLP}.

\subsection{Multistage blending problem}
\label{sec:blending}

\begin{figure}[ht]
\centering    \includegraphics[width=0.5\linewidth]{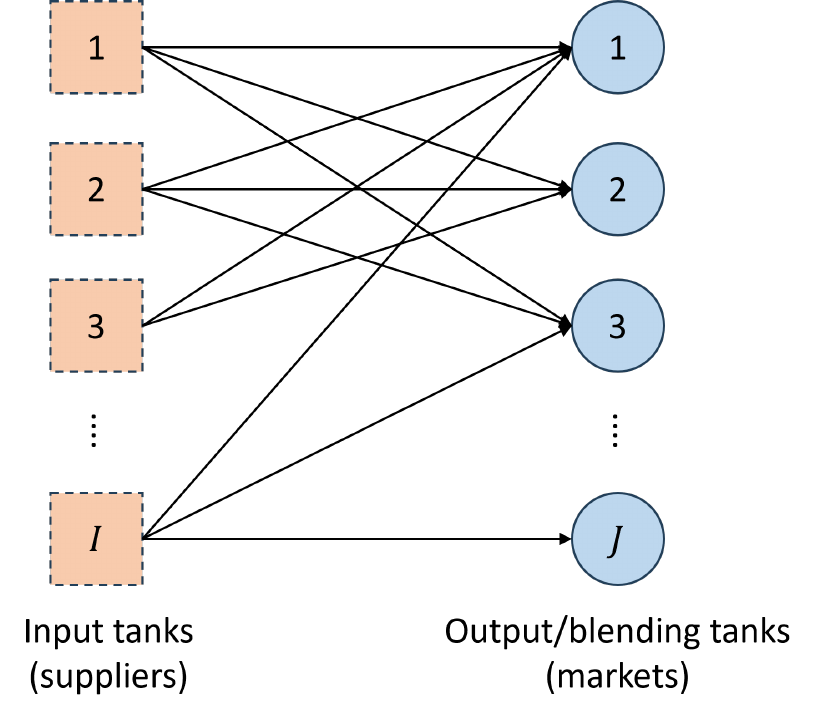}
    \caption{A network of input (suppliers) and output (markets) tanks. The specification of the product delivered to each market is determined by the blending of input streams from different suppliers.}
    \label{blending_network}
\end{figure}

Blending or pooling problems \citep{misener2009advances} appear in a variety of industrial processes, including oil refining, chemical manufacturing, food and beverage production, and the supply chains of various commodities. In this case study, we assume a general network (Figure \ref{blending_network}) with potential input tanks and a fixed number of output or blending tanks. The output tanks represent commodity markets and are assumed to offer different prices depending on the specifications of the supplied products. The input tanks are representative of suppliers and are installed as needed to meet demand. The objective is driven by the need to maximize overall revenue over a given planning horizon.

\subsubsection{Model formulation}

The deterministic model for the multiperiod blending problem is as follows:
\begin{align}
\maximize_{x,c,d,F} {}&\quad \sum\limits_{t \in \mathcal{T}}\Big[\sum_{i \in \mathcal{I}}\sum_{j \in \mathcal{J}}f_{jt}(c_{jt}) F_{ijt} - \sum\limits_{i \in \mathcal{I}}q_{it}x_{it} - \sum_{i \in \mathcal{I}}\sum_{j \in \mathcal{J}}r_{ij}F_{ijt} - \sum_{i \in \mathcal{I}}b_{it}\sum_{j \in \mathcal{J}}F_{ijt}\Big] \label{cs1_objective} &\\
\sto & \;\; \sum\limits_{i \in \mathcal{I}}\lambda_{i}F_{ijt} = c_{jt}d_{jt} \quad \forall \, j \in \mathcal{J}, t \in \mathcal{T} \label{cs1_blending}\\
& \;\; \sum\limits_{i \in \mathcal{I}}F_{ijt} = d_{jt} \quad \forall \, j \in \mathcal{J}, t \in \mathcal{T} \label{cs1_demand}\\
& \;\; \sum\limits_{j \in \mathcal{J}}F_{ijt} \leq C_{i} \sum\limits_{\tau=1}^{t}x_{i\tau} \quad \forall \, i \in \mathcal{I}, t \in \mathcal{T} \label{cs1_capacity1}\\
& \;\; \sum\limits_{\tau=1}^{t}x_{i\tau} \leq 1 \quad \forall \, i \in \mathcal{I}, t \in \mathcal{T} \label{cs1_capacity2}\\
& \;\; 0 \leq c_{jt} \leq 1 \quad \forall \, j \in \mathcal{J}, t \in \mathcal{T} \label{cs1_bound_c}\\
& \;\; d_{jt}^{\mathrm{min}} \leq d_{jt} \leq d_{jt}^{\mathrm{max}} \quad \forall \, j \in \mathcal{J}, t \in \mathcal{T} \label{cs1_bound_d}\\
& \;\; F_{ijt} \geq 0 \quad \forall \, i \in \mathcal{I}, j \in \mathcal{J}, t \in \mathcal{T} \label{cs1_bound_F}\\
& \;\; x_{it} \in \{0,1\} \quad \forall \, i \in \mathcal{I}, t \in \mathcal{T}, \label{cs1_bound_x}
\end{align}
where sets $\mathcal{I}$ and $\mathcal{J}$ represent the input and output tanks, whereas set $\mathcal{T}$ represents the time periods in the planning horizon. The binary variable $x_{it}$ equals 1 if input tank $i$ is installed in time period $t$. The product flow between input tank $i$ and output tank $j$ in time period $t$ is denoted by $F_{ijt}$. The product demand satisfied at output tank $j$ in time period $t$ is represented by $d_{jt}$ and is bounded by $d_{jt}^{\mathrm{min}}$ and $d_{jt}^{\mathrm{max}}$. The quality specification of the product that market $j$ receives (post blending) in time period $t$ is denoted by $c_{jt}$. The quality specification of the product manufactured at supplier $i$ is denoted by $\lambda_{i}$. The unit production cost at input tank $i$ in time period $t$ and the unit transportation cost from input tank $i$ to output tank $j$ are denoted by $b_{it}$ and $r_{ij}$, respectively. The cost of installing input tank $i$ of capacity $C_{i}$ in time period $t$ is denoted by $q_{it}$. The unit price offered by output tank $j$ for a product of specification $c_{jt}$ in time period $t$ is captured by $f_{jt}(c_{jt})$. For this case study, we assume that $f(c)$ is linear, i.e. of the form $mc+l$, where $l$ represents the unit cost for a product of specification $c=0$, and $m$ represents the factor by which the offered price increases (or decreases) with the specification. Constraints \eqref{cs1_blending} enforce the linear blending of the product supplied to each output tank. Constraints \eqref{cs1_demand} ensure demand satisfaction at each output tank. Constraints \eqref{cs1_capacity1} ensure that the net flow out of an input tank does not exceed its capacity. Constraints \eqref{cs1_capacity2} ensure that an input tank is installed only once during the planning horizon. Constraints \eqref{cs1_bound_c}-\eqref{cs1_bound_x} indicate the bounds on the decision variables. Lastly, according to \eqref{cs1_objective}, the objective is to maximize the overall revenue generated during the planning horizon. Assuming $f(c)$ is a linear function, the above model is a nonconvex MINLP due to the bilinear terms in the objective function and the blending constraints \eqref{cs1_blending}. The input tank installation decisions $x$ are the state variables that link the time periods.

We consider stochasticity in the minimum and maximum demands at each output tank, represented by the parameters $d_{jt}^{\mathrm{min}}$ and $d_{jt}^{max}$, respectively. To fit the decomposable structure of (MSSP), the resulting stochastic programming formulation requires the introduction of an auxiliary variable, which also enables the column sharing step to generate meaningful columns and is detailed in Appendix A. Appendix B contains the master problem and subproblem formulations for the column generation decomposition.

\subsubsection{Results and discussion}

We now generate instances of different sizes to assess and compare the performance of the proposed decomposition approach over solving a fullspace model. Specifically, the instances correspond to having (a) 5 input, 3 output, or (b) 12 input, 10 output tanks. Further, within each case, the number of time periods ($|\mathcal{T}|$) was varied from 3 to 9 in increments of 2, i.e. 3, 5, 7, and 9, which correspond to 8, 32, 128, and 512 scenarios ($|\mathcal{S}|$), respectively. The model size statistics are summarized in Table \ref{tab:blending_model_size_stats}. All instances were run on Minnesota Supercomputing Institute's AMD EPYC 7702 Linux cluster Mangi. A termination criterion of 0.01\% optimality gap with a 2-hour (7,200 s) time limit was set for each instance.

\begin{table}[ht]
\fontsize{12}{16}\selectfont
\setlength\tabcolsep{5pt}
  \centering
\begin{tabular}{clrrr}
\hline
\multicolumn{1}{l}{\textbf{($\bm{|\mathcal{I}|,|\mathcal{J}|}$)}} & \textbf{$\bm{|\mathcal{T}|/|\mathcal{S}|}$} & \textbf{\# of binary vars.} & \textbf{\# of continuous vars.} & \textbf{\# of constraints} \\ \hline
\multirow{4}{*}{(5,3)} & 3/8 & 35 & 294 & 597 \\
 & 5/32 & 155 & 1,302 & 2,661 \\
 & 7/128 & 635 & 5,334 & 10,917 \\
 & 9/512 & 2,555 & 21,462 & 43,941 \\ \hline
\multirow{4}{*}{(12,10)} & 3/8 & 84 & 1,960 & 2,844 \\
 & 5/32 & 372 & 8,680 & 12,636 \\
 & 7/128 & 1,524 & 35,560 & 51,804 \\
 & 9/512 & 6,132 & 143,080 & 208,476 \\ \hline
\end{tabular}%
    \caption{Model size statistics. Note that the number of constraints includes variable bounds.}
  \label{tab:blending_model_size_stats}%
\end{table}

Five instances were solved for each combination of the number of time periods ($|\mathcal{T}|$) and the numbers of input ($|\mathcal{I}|$) / output ($|\mathcal{J}|$) tanks. The distributions from which the model parameters were sampled are provided in Appendix E. The performance statistics are shown in Table \ref{tab:blending_result}. On the smaller instances with 3 and 5 time periods, the fullspace model performs well; on the larger instances with 7 and 9 time periods, however, its performance significantly deteriorates. Specifically, fullspace is unable to solve 35\% of larger instances with 7 and 9 time periods to optimality; in these cases, the final mean gap is between approximately 6\% and 13\%. CG outperforms the fullspace model in terms of the eventual optimality gap of the instances that are not solved to optimality. Even though the performance improvement may not seem substantial at first (6 problems were not solved to optimality here over 7 in the case of the fullspace model), the mean final gap is between 0.04\% and 0.69\%, which is considered optimal for many real-world scenarios. Furthermore, recall that we utilized a limited number of cores to parallelize the independent subproblems. As high-performance computing systems become more accessible and allow for node sharing, it is becoming increasingly possible to use a larger number of cores to \textit{perfectly} parallelize the subproblems, with all subproblems allocated to different cores that start solving at the same time (which may not always be possible due to shared resources and their availability). This can significantly increase the computational efficiency of parallelizable algorithms. As a result, we also provide an optimistic estimate of the solution time under perfectly parallelizable conditions, which can lead to significant speed gains. For example, for instances with 12 input and 10 output tanks that were solved to optimality, the solution time is estimated to improve by $\sim$ 23\% and 77\% in the 7 and 9 time period cases, respectively.

\begin{table}[ht]
\fontsize{12}{16}\selectfont
\setlength\tabcolsep{2pt}
  \centering
\resizebox{\columnwidth}{!}{%
\begin{tabular}{@{}clrrrrrrrrrrr@{}}
\toprule
\multicolumn{1}{l}{} &  & \multicolumn{3}{c}{\textbf{Fullspace}} & \multicolumn{4}{c}{\textbf{CG}} & \multicolumn{4}{c}{\textbf{CGCS}} \\ \cmidrule(lr){3-5} \cmidrule(lr){6-9} \cmidrule{10-13} 
\multicolumn{1}{l}{\textbf{($\bm{|\mathcal{I}|,|\mathcal{J}|}$)}} & \textbf{$\bm{|\mathcal{T}|/|\mathcal{S}|}$} & \multicolumn{1}{l}{\textbf{NS}} & \multicolumn{1}{l}{\textbf{\textoverline{gap} (\%)}} & \textbf{time (s)} & \multicolumn{1}{l}{\textbf{NS\textsuperscript{*}}} & \multicolumn{1}{l}{\textbf{\textoverline{gap} (\%)}} & \textbf{time (s)} & \textbf{time\textsuperscript{\#} (s)} & \multicolumn{1}{l}{\textbf{NS\textsuperscript{*}}} & \multicolumn{1}{l}{\textbf{\textoverline{gap} (\%)}} & \textbf{time (s)} & \textbf{time\textsuperscript{\#} (s)} \\ \midrule
\multirow{4}{*}{(5,3)} & 3/8 & 0 & - & 1 & 0 & - & 42 & 3 & 0 & - & 61 & 3 \\
 & 5/32 & 0 & - & 6 & 0 & - & 37 & 18 & 0 & - & 41 & 18 \\
 & 7/128 & 0 & - & 163 & 0 & - & 151 & 61 & 0 & - & 120 & 43 \\
 & 9/512 & 2 & 5.75 & 2,273 & 1 & 0.69 & 669 & 80 & 1 & 0.37 & 603 & 57 \\ \midrule
\multirow{4}{*}{(12,10)} & 3/8 & 0 & - & 8 & 0 & - & 199 & 156 & 0 & - & 145 & 88 \\
 & 5/32 & 0 & - & 148 & 0 & - & 684 & 640 & 0 & - & 372 & 347 \\
 & 7/128 & 3 & 12.87 & 83 & 1 & 0.04 & 2,733 & 2,104 & 1 & 0.04 & 1,258 & 995 \\
 & 9/512 & 2 & 8.88 & 1,239 & 4 & 0.19 & 5,428 & 1,241 & 0 & - & 5,368 & 933 \\ \bottomrule
\end{tabular}%
}
      \caption{Summary statistics highlighting the differences in performance of the fullspace model, column generation (CG), and column generation with column sharing (CGCS). For every combination of $(\mathcal{I},\mathcal{J})$ and $\mathcal{T}/\mathcal{S}$, 5 random instances were solved, and the average statistics are reported. ($|\mathcal{I}|$: number of input tanks, $|\mathcal{J}|$: number of output tanks, $|\mathcal{T}|$: number of time periods, $|\mathcal{S}|$: number of scenarios,
      NS: number of instances not solved to 0.01\% optimality gap in 7,200 s, NS\textsuperscript{*}: number of instances not solved to optimality in 7,200 s or column generation converged under 7,200 s but (RM-MSSP) did not provide a solution under 0.01\% optimality gap, \textoverline{gap}: average optimality gap for instances not solved to optimality, time: average solution time for instances solved to 0.01\% optimality gap,
      time\textsuperscript{\#}: average solution time for instances solved to 0.01\% optimality gap assuming perfect parallelization.)}
  \label{tab:blending_result}%
\end{table}

\begin{remark}
    In the perfectly parallelizable case, the solution time is estimated by assuming that the time spent in the pricing step in each iteration equals the time spent solving the slowest subproblem (+ the slowest CS problem in the CGCS case), simulating the situation in which all subproblems are solved on different cores simultaneously. It should be noted that this optimistic estimate of solution time only provides a rough idea of the degree of improvement that can be achieved and is not necessarily a lower bound, and it may vary considerably depending on the number of actual cores made available, shared memory, machine specifications, and method of parallelization (e.g., asynchronous solving of regular subproblems and the CS problems), among many other factors. 
\end{remark}

Next, we see additional benefits from incorporating the column sharing procedure into the column generation algorithm. In particular, more instances can be solved to optimality, and the solution time reduces considerably, especially for larger instances, as shown in Table \ref{tab:blending_result}. The first key observation is that for the 9 time period case with 12 input and 10 output tanks, CGCS solves all problems to optimality, whereas CG failed at closing the gap in 4 out of 5 instances. Second, CGCS leads to significantly lower solution times, especially for the 5 and 7 time period cases in the 12 input and 10 output cases. For the 5 time period case, the average solution time drops by $\sim$ 46\% from 684 s to 372 s. Similarly, for the 7 time period case, the solution time improves by $\sim$ 54\% from 2,733 s to 1,258 s. Again, similar to the observation in the case of CG, perfect parallelization can further enhance the performance of CGCS for both smaller as well as larger instances. For instances that CGCS solved to optimality, the estimated mean solution time improvement ranged from $\sim$ 56\% to 95\% and $\sim$ 7\% to 83\% in the 5 input/3 output and 12 input/10 output cases, respectively.

\begin{figure}[ht]
    \centering
    \includegraphics[width=0.8\linewidth]{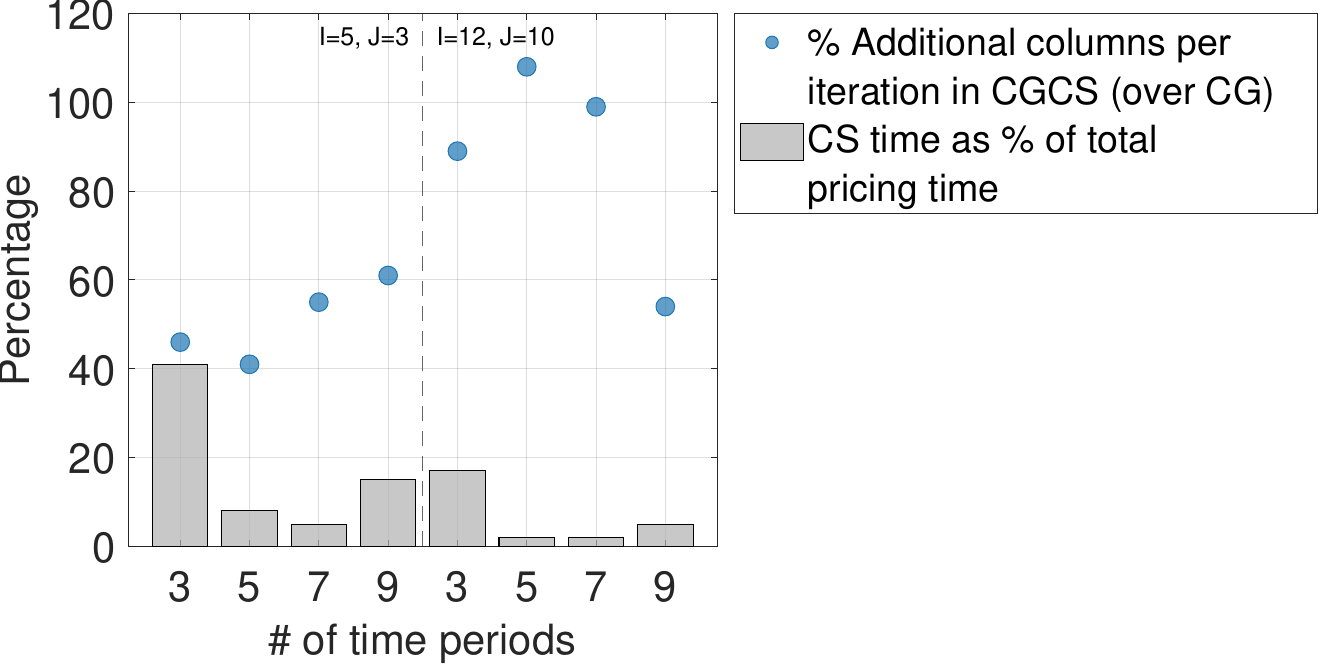}
    \caption{Illustrating the percentage of additional columns generated in CGCS versus CG, as well as the percentage of time spent in the CS step compared to the total pricing time (solving regular subproblems + CS).}
    \label{fig:blending_cols_time}
\end{figure}

\begin{figure}[ht]
    \centering
    \begin{subfigure}[b]{0.45\linewidth}
        \centering
        \includegraphics[width=\linewidth]{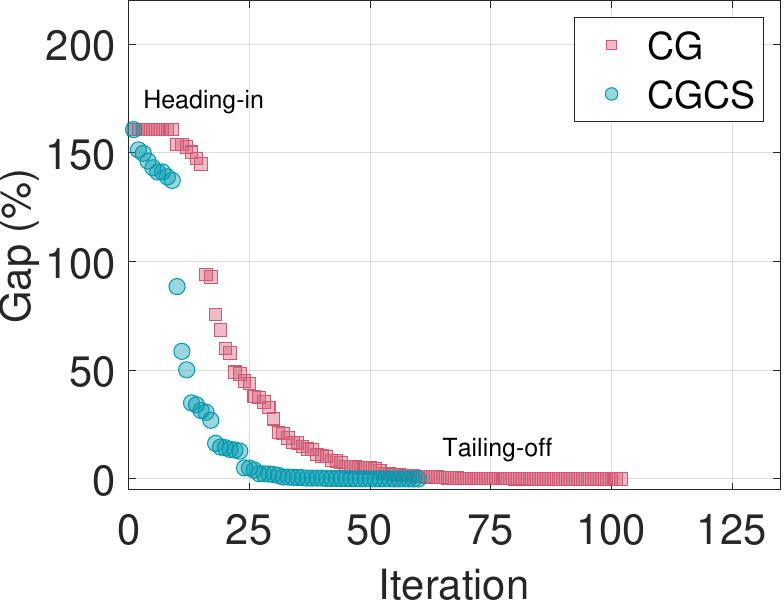}
        \captionsetup{justification=centering, margin={0cm,-1cm}}
        \caption{}
        \label{fig:blending_conv_plot}%
    \end{subfigure}
    \hfill
    \begin{subfigure}[b]{0.45\linewidth}
        \centering
        \includegraphics[width=\linewidth]{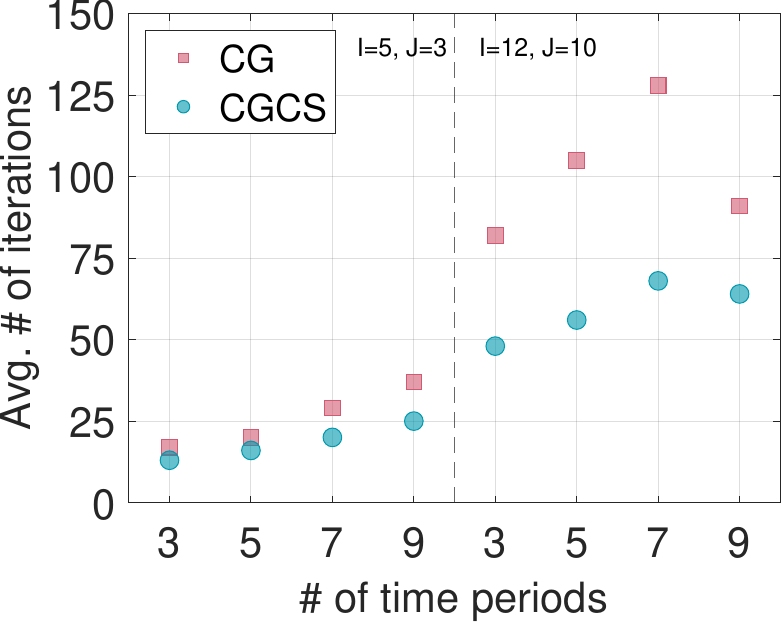}
        \captionsetup{justification=centering, margin={0cm,-1cm}}
        \caption{}
        \label{fig:blending_iter}
    \end{subfigure}
    \caption{(a) Convergence profile for an instance from 12 input/10 output with 7 time period case, and (b) average number of iterations for instances of different sizes.}
    \label{fig:comparison}
\end{figure}

Now we investigate why CGCS outperforms CG. We highlight two key correlated aspects here: (i) As expected, CGCS generates significantly more columns per iteration than CG. As shown in Figure \ref{fig:blending_cols_time}, the additional columns generated in CGCS over CG varied from 41\% to 108\%. Recall that these additional columns naturally satisfy the nonanticipativity constraints, increasing the likelihood of generating more feasible solutions and preventing the algorithm from stalling. The convergence profile in Figure \ref{fig:blending_conv_plot} for one of the instances from the 12 input/10 output with 7 time period case demonstrates this. In contrast to CG, CGCS eliminates the heading-in effect, results in faster convergence in intermediate iterations, and has a less pronounced tailing-off effect. This results in requiring fewer iterations to converge, as is highlighted in Figure \ref{fig:blending_iter}. (ii) Although CS is a step in addition to solving the regular subproblems, it does not necessarily slow down the algorithm. As shown in Figure \ref{fig:blending_cols_time}, in most cases, CS accounts for less than 20\% of the total pricing time (solving regular subproblems + CS), and in many cases, it can be as low as 2-5\%. Recall that in CS, a large number of integer variables are fixed (columns from sibling nodes), thus reducing its computational complexity. In conclusion, slight overhead in the pricing step, if any, is compensated by additional high-quality columns generated by CS in every iteration, effectively reducing the overall solution time.

\subsection{Mobile generator routing problem}
\label{sec:OPF}

Optimal power flow formulations that incorporate the modeling of power flow physics are frequently used to aid in the planning, scheduling, coordination, and control of the power grid \citep{Frank2012a}. One such application is to optimize the operation of distribution networks in a way that ensures power can be delivered to all customers. This task is made especially challenging by factors such as uncertainty in customer demand (load) forecasts and intermittent changes in customer behavior and power usage (e.g., due to HVAC usage, weeks-long or months-long industrial or agricultural activities, intermittent tourism, etc.). One option for addressing these challenges is the use of mobile generators (gensets) as non-wires alternatives for power delivery \citep{Riley2023}. These gensets are able to provide flexible backup generation, which helps address uncertainty in load forecasts, and can be relocated as power usage patterns change. However, the routing of these devices and the operation of the distribution network with these devices still yields a significant challenge, as discrete decisions (the location of gensets) need to be made while considering nonlinear power flow physics and uncertain operating conditions. 

\subsubsection{Model formulation}
\begin{figure}[H]
\centering    
\includegraphics[width=0.8\linewidth]{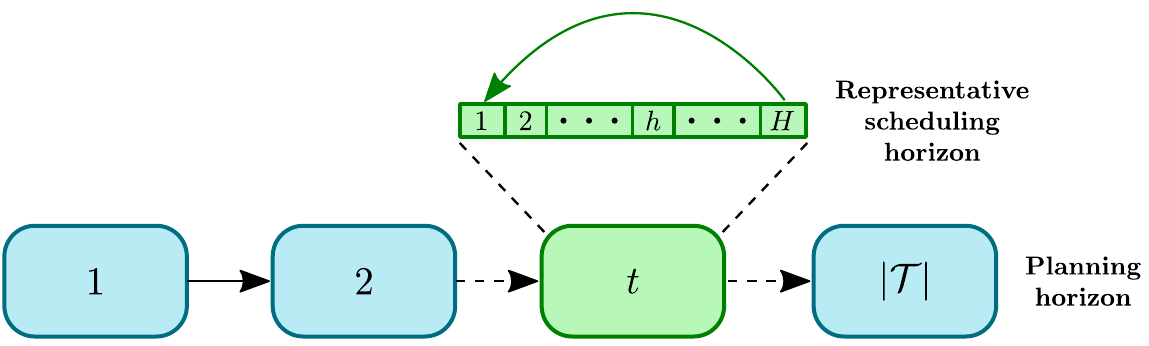}
    \caption{Multiscale time representation, which divides the planning horizon into a set of time periods, $\mathcal{T}$, with each time period $t \in \mathcal{T}$ having a representative scheduling horizon $\mathcal{H}$ of length $H$. In this case study, each time period is one week long, with a representative scheduling horizon of 24 hours.}
\label{fig:MSR}
\end{figure}

\begin{figure}[H]
\centering    
\includegraphics[width=0.8\linewidth]{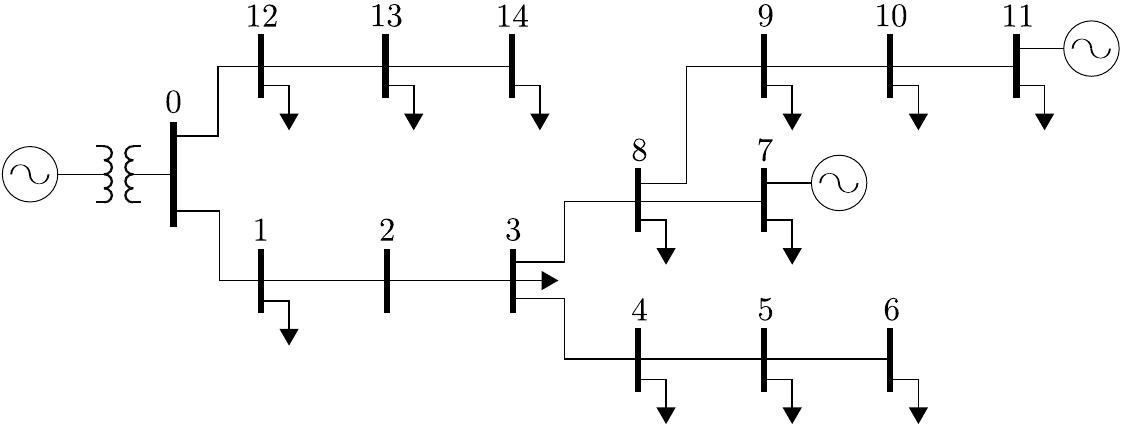}
    \caption{15-bus distribution system used in the case study.}
    \label{fig:power_network}
\end{figure}

In this case study, we solve a multistage stochastic program that minimizes the expected operating cost of a radial power distribution network by deciding where to locate gensets over time to respond to fluctuating loads at different buses in the network. In addition, we use a multiscale time representation (Figure \ref{fig:MSR}) to efficiently model hourly operational decisions over a given planning horizon. Furthermore, the multiscale time representation enables the efficient modeling of both hourly and longer-term, intermittent variations in load. Figure \ref{fig:power_network} depicts the schematic of a modified 15-bus distribution system \citep{Papavasiliou2018} that is considered in this case study. The distribution system includes two distributed stationary sources of generation located at Buses 7 and 11. Additionally, one or two mobile gensets (depending on problem instance) are available in the system.

\begin{figure}[H]
\centering    
\includegraphics[width=0.8\linewidth]{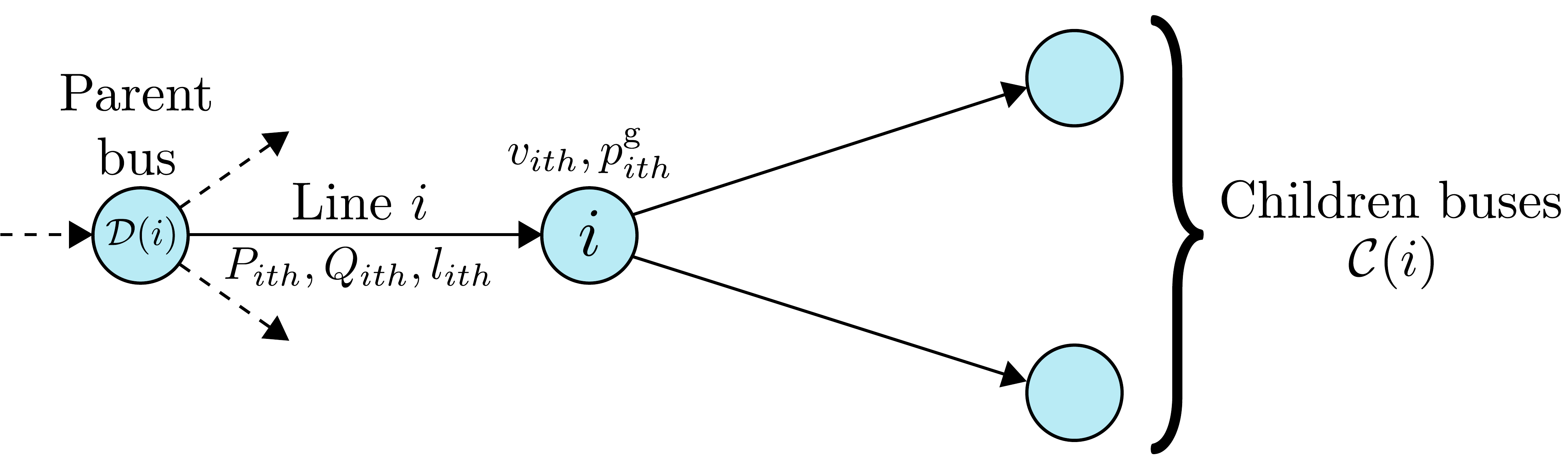}
    \caption{Directed graph representation of power distribution network, with bus $i$, its parent and children buses, line $i$, and essential notation highlighted.}
    \label{fig:power_network_notation}
\end{figure}

The deterministic multiperiod mobile generator routing problem employs the conic relaxation of the relaxed branch flow model \citep{farivar2013branch} and is formulated as follows:
\begin{align}
\minimize {}&\quad \sum_{t \in \mathcal{T}\backslash\{1\}}\sum_{m \in \mathcal{M}}C^\mathrm{t}_{mt}z_{mt}  + \label{cs2_objective}\\
& \qquad \qquad \sum\limits_{t \in \mathcal{T}} \alpha_{t} \sum_{i \in \mathcal{B}} \sum_{h \in \mathcal{H}} \bigg[C_{ith}^{\mathrm{p}}p_{ith}^{\mathrm{g}}+ C^{\mathrm{VoLL}}_i (1-\delta_{ith}) p^{\mathrm{d}}_{ith} + \sum_{m \in \mathcal{M}}{C^{\mathrm{a}}_{t}p_{mith}^{\mathrm{a}}}\bigg]  \nonumber\\
\sto & \quad \sum_{i \in \mathcal{B}}y_{mit}=1 \quad \forall \, m \in \mathcal{M}, t \in \mathcal{T} \label{cs2_sumy}\\
& \quad y_{mit}-y_{m,i,t-1} \leq z_{mt} \quad \forall \, m \in \mathcal{M}, i \in \mathcal{B}, t \in \mathcal{T}\backslash\{1\} \label{cs2_z1}\\
& \quad z_{mt} \leq 2 - y_{mit} - y_{m,i,t-1} \quad \forall \, m \in \mathcal{M}, i \in \mathcal{B}, t \in \mathcal{T}\backslash\{1\} \label{cs2_z2}\\
& \quad p_{mith}^{\mathrm{a}} \leq \overline{p}_{m}^{\mathrm{a}} y_{mit} \quad \forall \, m \in \mathcal{M}, i \in \mathcal{B}, t \in \mathcal{T}, h \in \mathcal{H} \label{cs2_gen_install2} \\
& \quad  -p_{mith}^{\mathrm{a}} \leq q_{mith}^{\mathrm{a}} \leq p_{mith}^{\mathrm{a}} \quad \forall \, m \in \mathcal{M}, i \in \mathcal{B}, t \in \mathcal{T}, h \in \mathcal{H} \label{cs2_gen_install3} \\
& \quad p_{ith}^{\mathrm{g}} - p_{ith}^{\mathrm{d}} \delta_{ith} + \sum_{m \in \mathcal{M}}p_{mith}^{\mathrm{a}} = \nonumber \\ 
& \quad \qquad \qquad \sum_{j \in \mathcal{C}(i)} P_{jth} - \sum_{j \in \mathcal{D}(i)}{(P_{ith}-r_{i}l_{ith})} + g_{i}v_{ith} \quad \forall \, i \in \mathcal{B}, t \in \mathcal{T}, h \in \mathcal{H} \label{cs2_powerflow_1} \\
& \quad q_{ith}^{\mathrm{g}} - q_{ith}^{\mathrm{d}} \delta_{ith} + \sum_{m \in \mathcal{M}}q_{mith}^{\mathrm{a}} = \nonumber \\
& \quad \qquad \qquad \sum_{j \in \mathcal{C}(i)} Q_{jth} - \sum_{j \in \mathcal{D}(i)}{(Q_{ith}-x_{i}l_{ith})} + b_{i}v_{ith} \quad \forall \, i \in \mathcal{B}, t \in \mathcal{T}, h \in \mathcal{H} \label{cs2_powerflow_2}\\
& \quad v_{ith} = v_{d(i),t,h} - 2(r_{i}P_{ith} + x_{i}Q_{ith}) + (r_{i}^{2} + x_{i}^{2})l_{ith} \quad \forall \, i \in \mathcal{L}, t \in \mathcal{T}, h \in \mathcal{H} \label{cs2_powerflow_3} \\
& \quad l_{ith}v_{d(i),t,h} \geq P_{ith}^2 + Q_{ith}^2 \quad \forall \, i \in \mathcal{L}, t \in \mathcal{T}, h \in \mathcal{H} \label{cs2_powerflow_4} \\
& \quad P_{ith}^2 + Q_{ith}^2 \leq A_{i}^2 \quad \forall \, i \in \mathcal{L}, t \in \mathcal{T}, h \in \mathcal{H} \label{cs2_powerflow_5} \\
& \quad (P_{ith}-r_{i}l_{ith})^2 + (Q_{ith}-x_{i}l_{ith})^2 \leq A_{i}^2 \quad \forall \, i \in \mathcal{L}, t \in \mathcal{T}, h \in \mathcal{H} \label{cs2_powerflow_6}\\
& \quad 0 \leq p_{ith}^{\mathrm{g}} \leq \overline{p}_{i}^{\mathrm{g}} \quad \forall \, i \in \mathcal{B}, t \in \mathcal{T}, h \in \mathcal{H} \label{cs2_bounds_1}\\
& \quad 0 \leq q_{ith}^{\mathrm{g}} \leq \overline{q}_{i}^{\mathrm{g}} \quad \forall \, i \in \mathcal{B}, t \in \mathcal{T}, h \in \mathcal{H} \label{cs2_bounds_2}\\
& \quad  0 \leq \delta_{ith} \leq 1 \quad \forall \, i \in \mathcal{B}, t \in \mathcal{T}, h \in \mathcal{H} \label{delta_bounds}\\
& \quad \underline{v}_{i} \leq v_{ith} \leq \overline{v}_{i} \quad \forall \, i \in \mathcal{B}, t \in \mathcal{T}, h \in \mathcal{H} \label{cs2_bounds_3}\\
& \quad l_{ith} \geq 0 \quad \forall \, i \in \mathcal{L}, t \in \mathcal{T}, h \in \mathcal{H} \label{cs2_bounds_4}\\
& \quad  p_{mith}^{\mathrm{a}} \geq 0 \quad \forall \, m \in \mathcal{M}, i \in \mathcal{B}, t \in \mathcal{T}, h \in \mathcal{H} \label{cs2_bounds_5}\\
& \quad y_{mit} \in \{0,1\} \quad \forall \, m \in \mathcal{M}, i \in \mathcal{B}, t \in \mathcal{T} \label{cs2_bounds_6} \\
& \quad z_{mt} \in \{0,1\} \quad \forall \, m \in \mathcal{M}, t \in \mathcal{T}\backslash\{1\} \label{cs2_bounds_7} \\
\operatorname{where} & \quad p_{ith}^{\mathrm{d}} = p^{\mathrm{d,base}}_i \beta_{it} e_{th} \quad \forall \, i \in \mathcal{B}, t \in \mathcal{T}, h \in \mathcal{H} \label{pd_definition} \\
& \quad q_{ith}^{\mathrm{d}} = q^{\mathrm{d,base}}_i \beta_{it} e_{th} \quad \forall \, i \in \mathcal{B}, t \in \mathcal{T}, h \in \mathcal{H}, \label{qd_definition}
\end{align}

\noindent where $\mathcal{M}$, $\mathcal{B}$, $\mathcal{L}$ denote the sets of mobile generators, buses, and power lines in the network, respectively. The binary variable $y_{mit}$ equals 1 if genset $m$ is located at bus $i$ in time period $t$, and the binary variable $z_{mt}$ equals 1 if genset $m$ is relocated at the beginning of time period $t$. Here, the decision variables pertaining to the power system operation are as follows: $p_{mith}^{\mathrm{a}}$ is the active power generated by genset $m$, $q_{mith}^{\mathrm{a}}$ is the reactive power generated by the same genset, $p_{ith}^{\mathrm{g}}$ and $q_{ith}^{\mathrm{g}}$ are the active and reactive power supplied by stationary generating units at bus $i$, $\delta_{ith}$ is the fraction of power demand supplied at bus $i$,  $P_{ith}$ and $Q_{ith}$ are the active and reactive power flows on line $i$, $l_{ith}$ is the squared current magnitude on line $i$, and $v_{ith}$ is the squared voltage magnitude at bus $i$. Additionally, $\mathcal{C}(i)$ denotes the set of children nodes of bus $i$ in the directed graph representation of the network, while $\mathcal{D}(i)$ denotes the set of parent nodes of bus $i$. $\mathcal{D}(i)$ is a singleton containing $d(i)$, the one parent node of bus $i$, for all buses but Bus 0, which has no parent node. Figure  \ref{fig:power_network_notation} depicts the graph representation of the radial network and highlights essential notation.

The cost associated with relocating genset $m$ at time period $t$ is $C^\mathrm{t}_{mt}$. The variable generation cost of stationary generation units at bus $i$ is denoted by $C_{ith}^{\mathrm{p}}$. Similarly, $C^{\mathrm{VoLL}}_i$ is the penalty cost associated with load shedding at bus $i$, and $C^{\mathrm{a}}_{t}$ is the variable generation cost of mobile gensets. Constraints \eqref{cs2_sumy} ensure that a genset is located at exactly one bus in each time period $t$. Constraints \eqref{cs2_z1} and \eqref{cs2_z2} ensure that $z_{mt}$ track the relocation of the gensets. Additionally, the initial location of gensets in each scheduling horizon is not constrained. Constraints \eqref{cs2_gen_install2} and \eqref{cs2_gen_install3} ensure that gensets can only inject active and reactive power at the buses they are located at. Constraints \eqref{cs2_powerflow_1}-\eqref{cs2_powerflow_4} represent the power flows in the system using the aforementioned conic relaxation of the relaxed branch flow model \citep{farivar2013branch}, which accurately models single-phase AC power flows in radial networks \citep{Gan2015}. Constraints \eqref{cs2_powerflow_5} and \eqref{cs2_powerflow_6} impose thermal limits on the apparent power flows over the lines. In the constraints modeling the power flow, $r_i$ is the resistance of line $i$, $g_i$ is the shunt conductance at bus $i$, $x_i$ is the reactance of line $i$, $b_i$ is the negative of the shunt susceptance at bus $i$, and $A_i$ is the apparent power flow limit of line $i$.  Equations (\ref{pd_definition}) and (\ref{qd_definition}) show that the active and reactive power demand parameters, $p_{ith}^{\mathrm{d}}$ and $ q_{ith}^{\mathrm{d}}$, are calculated using products of several factors. Here, $p^{\mathrm{d,base}}_i$ and $q^{\mathrm{d,base}}_i$ represent the peak daily active and reactive power demand for a typical or average day, $\beta_{it}$ is a factor that accounts for the variability in power demands at bus $i$ during period $t$ due to intermittent effects, and $e_{th}$ is a factor accounting for the effect of the time of day on the power demand. The objective function \eqref{cs2_objective} minimizes the operating costs incurred during the planning horizon.  The sum of operating costs over representative scheduling horizon $\mathcal{H}$ is multiplied by the factor $\alpha_t$ to convert it to the total operating costs incurred during time period $t$. Note that the state variables are the binary variables $y$ that indicate the location of the gensets, and the nonlinearity in the model is due to the power flow constraints.

We model stochasticity in power demand at each bus using the load variation factor $\beta_{it}$. The resulting stochastic programming formulation, as well as the auxiliary variable introduced to enable meaningful column sharing, are detailed in Appendix C. Appendix D contains the master problem and subproblem formulations for the column generation decomposition. Note that the given MINLP is jointly convex in all, continuous and binary, variables; as such, there is a larger set of existing decomposition methods that can be applied to solve the problem. For example, it is amenable to the standard generalized Benders decomposition algorithm if its master problem is designed to include all discrete variables; its computational performance typically depends on the number of discrete variables, which increases with the number of scenarios and can be particularly large if many of the stage variables are integer.

\subsubsection{Results and discussion}
To evaluate the performance of the fullspace model and the column generation algorithm, we consider cases with one and two gensets. For each of these cases, we vary the number of time periods ($|\mathcal{T}|$), and thus the number of stages and scenarios ($|\mathcal{S}|$), generating problems of various sizes (Table \ref{tab:OPF_model_size_stats}). To initiate the CG/CGCS algorithm, we generate initial columns by solving an expected value problem. We use Gurobi's solution pool functionality, with \textit{PoolSearchMode} set to 1, to collect up to 50 diverse initial solutions. It should be noted that 50 different solutions do not necessarily imply 50 different columns for each node $n \in \mathcal{N}$. This is because two different solutions can have the same columns for the majority of the nodes in the scenario tree, resulting in significantly fewer initial feasible columns for each node $n$. Moreover, at each iteration of CG and CGCS, we use Gurobi's solution pool functionality (\textit{PoolSearchMode} set to 2)  to price out up to 5 columns with the lowest reduced costs, but only those with a negative reduced cost are included. In the CGCS method, all priced out columns with a negative reduced cost are subjected to column sharing. All models were solved on Minnesota Supercomputing Institute's AMD EPYC 7763 Linux cluster, Agate. A termination criterion of 0.1\% optimality gap with a time limit of 3 hours (10,800 s) was set for each instance.

\begin{table}[ht]
\fontsize{12}{16}\selectfont
\setlength\tabcolsep{5pt}
  \centering
\begin{tabular}{clrrr}
\hline
\multicolumn{1}{l} {$\bm{|\mathcal{M}|}$} & $\bm{|\mathcal{T}|/|\mathcal{S}|}$ & \textbf{\# of binary vars.} & \textbf{\# of continuous vars.} & \textbf{\# of constraints} \\ \hline
\multirow{4}{*}{1} & 2/9 & 63 & 38,016 & 80,797 \\
 & 3/36 & 207 & 152,064 & 323,140 \\
 & 4/81 & 639 & 380,160 & 808,249 \\
 & 5/128 & 1,871 & 772,992 & 1,645,148 \\ \hline
\multirow{4}{*}{2} & 2/9 & 126 & 46,656 & 98,234 \\
 & 3/36 & 414 & 186,624 & 392,840 \\
 & 4/81 & 1,278 & 466,560 & 982,898 \\
 & 5/128 & 3,742 & 948,672 & 2,001,976 \\ \hline
\end{tabular}%
    \caption{Model size statistics. Note that the number of constraints includes variable bounds.}
    \label{tab:OPF_model_size_stats}%
\end{table}

\begin{table}[ht]
\fontsize{12}{16}\selectfont
\setlength\tabcolsep{2pt}
  \centering
\resizebox{\columnwidth}{!}{%
\begin{threeparttable}
\begin{tabular}{@{}clrrrrrrrrrrr@{}}
\toprule
\multicolumn{1}{l}{} &  & \multicolumn{3}{c}{\textbf{Fullspace}} & \multicolumn{4}{c}{\textbf{CG}} & \multicolumn{4}{c}{\textbf{CGCS}} \\ 
\cmidrule(lr){3-5} \cmidrule(lr){6-9} \cmidrule{10-13} 
\multicolumn{1}{l}{$\bm{|\mathcal{M}|}$} & \textbf{$\bm{|\mathcal{T}|/|\mathcal{S}|}$} & \multicolumn{1}{l}{\textbf{NS}} & \multicolumn{1}{l}{\textbf{\textoverline{gap} (\%)}} & \multicolumn{1}{l}{\textbf{time (s)}} & \multicolumn{1}{l}{\textbf{NS}} & \multicolumn{1}{l}{\textbf{\textoverline{gap} (\%)}} & \textbf{time (s)} & \textbf{time\textsuperscript{\#} (s)} & \multicolumn{1}{l}{\textbf{NS}} & \multicolumn{1}{l}{\textbf{\textoverline{gap} (\%)}} & \textbf{time (s)} & \textbf{time\textsuperscript{\#} (s)} \\ \hline
\multirow{4}{*}{1} & 2/9 & 0 & - & \multicolumn{1}{r}{87} & 0 & - & 107 & 80 & 0 & - & 82 & 58 \\
 & 3/36 & 3 & 0.25 & \multicolumn{1}{r}{6,632} & 0 & - & \multicolumn{1}{r}{2,287} & \multicolumn{1}{r}{2,255} & 0 & - & \multicolumn{1}{r}{930} & \multicolumn{1}{r}{894} \\
 & 4/81 & 5 & 1.57 & - & 1 & 0.19 & 4,902 & 4,093 & 0 & - & 3,157 & 2,829 \\
 & 5/128 & \multicolumn{1}{r}{5} & \multicolumn{1}{r}{2.16} & - & \multicolumn{1}{r}{1} & 0.14 & 10,655 & 5,476 & \multicolumn{1}{r}{0} & - & 7,332 & 3,408 \\ \hline
\multirow{4}{*}{2} & 2/9 & 0 & - & 1,777 & 0 & - & 722 & 686 & 0 & - & 431 & 387 \\
 & 3/36 & 5 & 5.27 & - & 0 & - & 6,501 & 6,469 & 0 & - & 3,193 & 3,163 \\
 & 4/81 & 5* & - & - & \multicolumn{1}{r}{5} & 3.00 & \multicolumn{1}{r}{-} & \multicolumn{1}{r}{-} & \multicolumn{1}{r}{3} & 0.56 & 9,612 & 7,569 \\
 & 5/128 & 5* & - & - & \multicolumn{1}{r}{5} & 13.26 & \multicolumn{1}{r}{-} & \multicolumn{1}{r}{-} & 5 & 2.27 & \multicolumn{1}{r}{-} & \multicolumn{1}{r}{-} \\ \hline
\end{tabular}%
\vspace{0.1em}
\begin{tablenotes}
\item [\raisebox{-1ex}{\scalebox{1.5}{*}}] \hspace{-0.4em} \small{No feasible solution found in the allotted time limit}
\end{tablenotes}
\end{threeparttable}
}
      \caption{Summary statistics highlighting the differences in performance of the fullspace model, column generation (CG), and column generation with column sharing (CGCS). For every combination of $\mathcal{M}$ and $\mathcal{T}/\mathcal{S}$, 5 random instances were solved, and the average statistics are reported. ($|\mathcal{M}|$: number of gensets, $|\mathcal{T}|$: number of time periods, $|\mathcal{S}|$: number of scenarios,
      NS: number of instances not solved to 0.1\% optimality gap in 10,800 s, \textoverline{gap}: average optimality gap for instances not solved to optimality in 10,800 s, time: average solution time for instances solved to 0.1\% optimality gap,
      time\textsuperscript{\#}: average solution time for instances solved to 0.1\% optimality gap assuming perfect parallelization.)}
  \label{tab:OPF_result}%
\end{table}

To report average statistics, we generate five random instances for each combination of gensets ($|\mathcal{M}|$) and time periods ($|\mathcal{T}|$). The static bus and line related parameters, along with parameters sampled from different distributions, can be found in Appendix E. All performance statistics are summarized in Table \ref{tab:OPF_result}. For the 1 genset case: (i) With 2 time periods, fullspace and CG show similar performance. (ii) For the larger cases with 3, 4, and 5 time periods, CG clearly outperforms fullspace by solving 87\% of the instances to optimality, whereas fullspace could only converge for 13\% of the instances. For the 2 gensets case: (i) With 2 time periods, although both fullspace and CG could close the gap for all instances, CG converges in less than half the time taken by the fullspace model. (iii) With 3 time periods, CG could solve all instances to optimality, while fullspace failed to converge for any, resulting in a mean gap of 5.27\%. (iv) For the 4 and 5 time period cases, fullspace model fails to find any feasible solution in the 3-hour time limit, while CG could find solutions with mean gap of 3\% and 13\% for the two cases, respectively. Additionally, as was also the observation in the first case study, under conditions of perfect parallelization, we could further improve the CG solution time. For example, in the 1 genset case with 4 and 5 time periods, the solution time estimate suggests improvements of approximately 17\% and 49\%, respectively.

Next, we observe that one can achieve significant additional performance gains when using CGCS over CG. As shown in Table \ref{tab:OPF_result}, for the 1 genset case, CGCS solves all instances to optimality with significantly shorter solution times. For example, for 1 genset with 3 time periods, CGCS reduces the mean solution time to 930 s, a reduction of $\sim$ 59\% compared to CG's 2,287 s. In the 2 gensets case with 2 and 3 time periods, CGCS reduces the solution times by $\sim$ 40\% and 51\%, respectively, compared to regular CG. In the 2 gensets case with 4 time periods, CGCS solves a larger number of instances to optimality. Additionally, the instances that are not solved to optimality have a mean gap of only 0.56\%, which is acceptable for most practical purposes. Lastly, for the largest set of instances (2 gensets and 5 time periods), CGCS could solve all instances to a reasonably acceptable gap of 2.27\%, which is a significant improvement over CG's mean gap of 13.26\%.

\begin{figure}[ht]
    \centering
    \includegraphics[width=0.8\linewidth]{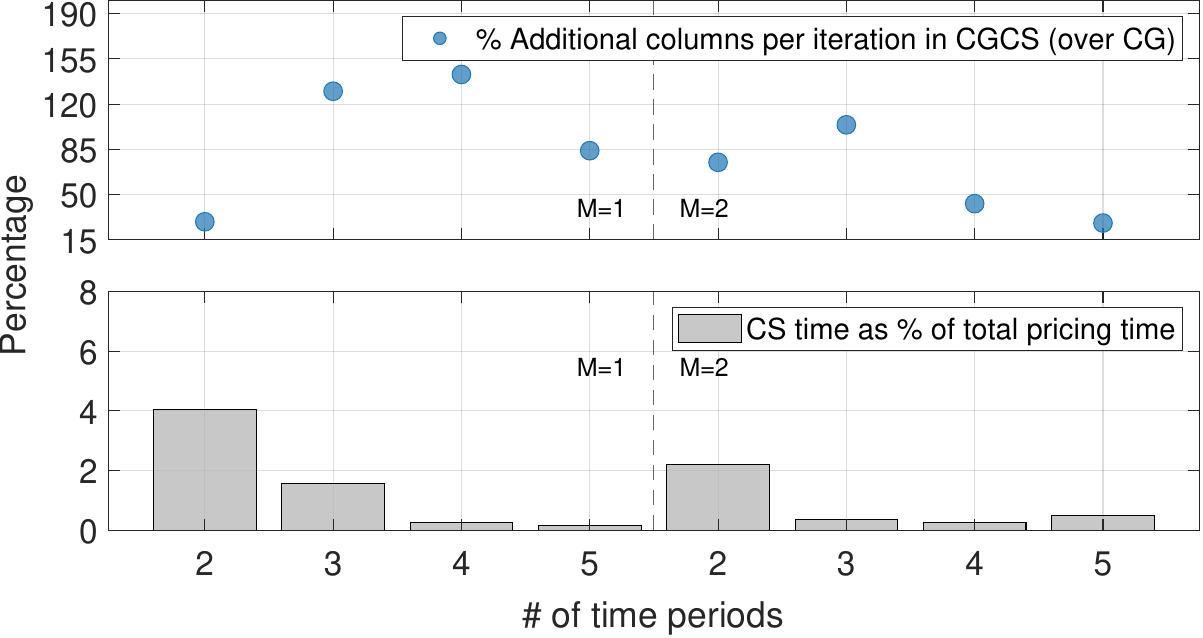}
    \caption{Illustrating the percentage of additional columns generated in CGCS versus CG, as well as the percentage of time spent in the CS step compared to the total pricing time (solving regular subproblems + CS).}
    \label{fig:opf_cols_time}
\end{figure}

\begin{figure}[ht]
    \centering
    \begin{subfigure}[b]{0.45\linewidth}
        \centering
        \includegraphics[width=\linewidth]{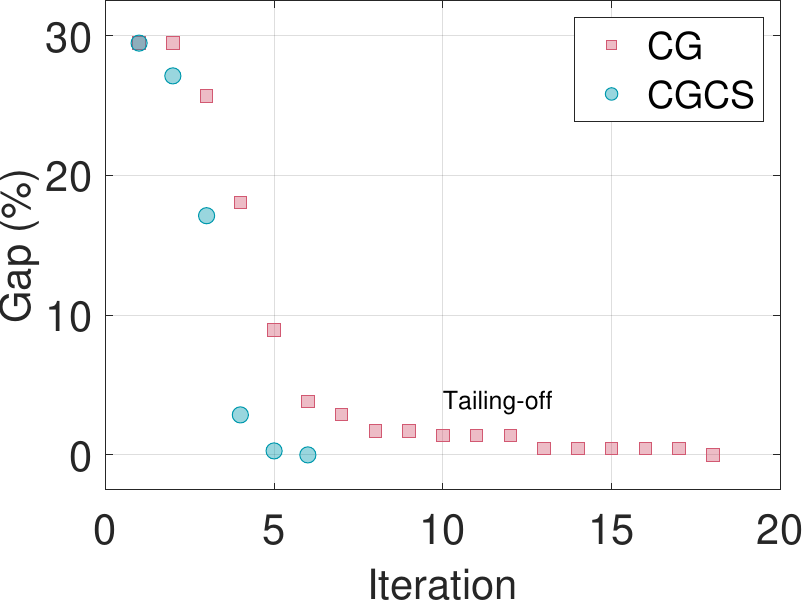}
        \captionsetup{justification=centering, margin={0cm,-1cm}}
        \caption{}
        \label{fig:opf_conv_plot-12}%
    \end{subfigure}
    \hfill
    \begin{subfigure}[b]{0.45\linewidth}
        \centering
        \includegraphics[width=\linewidth]{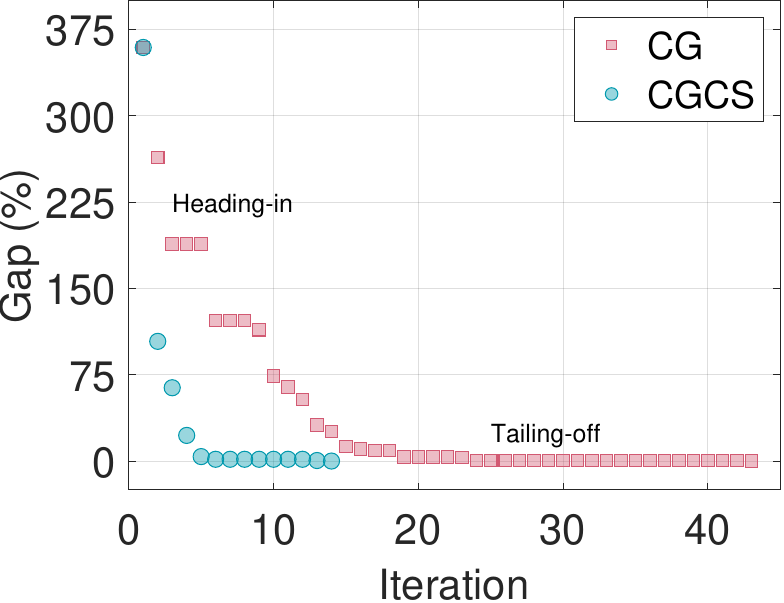}
        \captionsetup{justification=centering, margin={0cm,-1cm}}
        \caption{}
        \label{fig:opf_conv_plot-27}
    \end{subfigure}
    \vskip\floatsep 
    \begin{subfigure}[b]{0.45\linewidth}
        \centering
        \includegraphics[width=\linewidth]{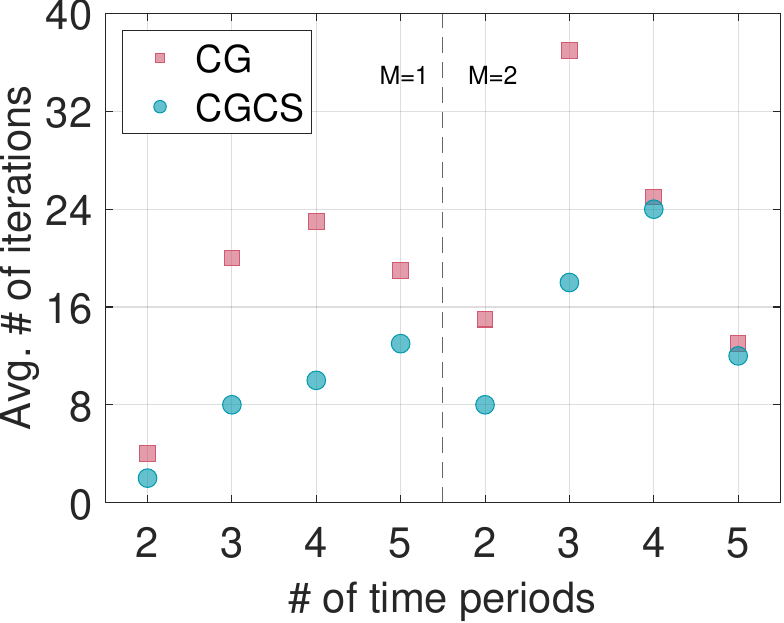}
        \captionsetup{justification=centering, margin={0cm,-1cm}}
        \caption{}
        \label{fig:opf_iter}
    \end{subfigure}
    \caption{Convergence profiles for (a) an instance from the 1 genset/4 time periods case, (b) an instance from the 2 gensets/3 time periods case. (c) The average number of iterations for different-sized instances until convergence or reaching the time limit.}
    \label{fig:opf-comparison}
\end{figure}

Now, like in the multistage blending case study, we present some arguments that highlight some of the reasons for CGCS' superior performance over CG. First, owing to the column sharing step in CGCS, we are able to generate significantly more columns per iteration as is shown in Figure \ref{fig:opf_cols_time}. For the instances we consider, CGCS generates between 28\% to 143\% additional columns per iteration. Since column sharing ensures nonanticipativity constraints are satisfied, these columns are more likely to form feasible solutions, decreasing the likelihood of algorithm stalling at various points. This is verified from the convergence plots shown in Figures \ref{fig:opf_conv_plot-12} and \ref{fig:opf_conv_plot-27} for a 1 genset/4 time periods case and a 2 gensets/3 time periods case, respectively. Clearly, CGCS exhibits faster convergence compared to CG. Specifically, in addition to largely mitigating the tailing-off effect in both cases, CGCS also eliminates the significant heading-in effect observed in the convergence under CG for the instance with 2 gensets (see Figure \ref{fig:opf_conv_plot-27}), ultimately requiring far fewer iterations to converge to an optimal solution. Figure \ref{fig:opf_iter} compares the mean number of iterations each algorithm underwent before converging or reaching the time limit. The only ambiguity is in the case of 2 gensets with 4 and 5 time periods. Here, we see both CG and CGCS undergo almost the same number of iterations; however, recall from Table \ref{tab:OPF_result} that CGCS converged to a much smaller gap in these many iterations. In the absence of a time limit, CGCS would have required fewer iterations than CG to reach the optimal solution. While a lower number of iterations does not always imply faster convergence, it is the case here. Figure \ref{fig:opf_cols_time} shows that the time spent on column sharing is often less than 0.5\% of the total pricing time for larger instances. This insignificant overhead in each iteration is offset by a significant reduction in total iterations, effectively reducing the overall solution time. 

\section{Conclusions}
\label{sec:conclusions}
In this work, we developed a column generation algorithm for solving generally nonconvex multistage stochastic MINLPs with discrete state variables by leveraging the specific structure of this class of problems. In addition, we applied a column sharing strategy to address known convergence issues in column generation.
We demonstrated the effectiveness of the proposed method by solving various instances of multistage blending and mobile generator routing problems. In both case studies, we observed that column generation consistently solves large instances to substantially lower optimality gaps than the fullspace model. Notably, for the large instances of the mobile generator routing problem, column generation yielded solutions with reasonable gaps (often near-optimal), whereas the fullspace model often failed to produce any feasible solution within the allotted time limit. Furthermore, integrating the column sharing procedure into the regular column generation algorithm resulted in significant performance improvements, both in terms of the number of optimally solved instances and convergence speed. Our detailed analysis reveals that the main reasons for this improvement are that column sharing enables the generation of additional columns in every iteration, all of which satisfy the nonanticipativity constraints, and that the column sharing procedure itself is computationally inexpensive. These factors work together to mitigate the convergence issues in column generation, resulting in more optimally solved instances and reduced solution times.

\section*{Acknowledgements}
The authors gratefully acknowledge the financial support from the State of Minnesota through an appropriation from the Renewable Development Account to the University of Minnesota as well as from the University of Minnesota through a MnDRIVE Research Computing Seed Grant. Angela Flores-Quiroz further acknowledges financial support from Instituto Sistemas Complejos de Ingenieria (Grant ANID PIA/PUENTE AFB230002). Computational resources were provided by the Minnesota Supercomputing Institute at the University of Minnesota.

\bibliographystyle{revunsrtnat}
\bibliography{library}

\newpage
\section*{Appendix}

\subsection*{A. Stochastic formulation for the multistage blending case study}
Here we present the stochastic programming formulation that models uncertainty in the minimum ($d_{jt}^{\mathrm{min}}$) and maximum ($d_{jt}^{\mathrm{max}}$) demands at each output tank.
The decision sequence in a time period is as follows: input tank installation decisions are made at the beginning of the time period, followed by the realization of uncertainty in the bounds on demands at each output tank, and finally, the product flow decisions are made. The deterministic formulation in Section \ref{sec:blending} is extended to the stochastic programming formulation as follows: 
\begin{align} 
\maximize_{x,c,d,F} {}&\quad \sum\limits_{n \in \mathcal{N}\backslash\{1\}}p_{n}\Big[\sum_{i \in \mathcal{I}}\sum_{j \in \mathcal{J}}f_{j,t(n)}(c_{jn}) F_{ijn} - \sum_{i \in \mathcal{I}}\sum_{j \in \mathcal{J}}r_{ij}F_{ijn} -  & \\
& \qquad \qquad \qquad \qquad \sum_{i \in \mathcal{I}}b_{i,t(n)}\sum_{j \in \mathcal{J}}F_{ijn}\Big] - \sum\limits_{n \in \mathcal{N}\backslash\mathcal{N}^{\mathrm{L}}}p_{n}\sum\limits_{i \in \mathcal{I}}q_{i,t(n)}x_{in} \nonumber \\
\sto & \;\; \sum\limits_{i \in \mathcal{I}}\lambda_{i}F_{ijn} = c_{jn}d_{jn} \quad \forall \, j \in \mathcal{J}, n \in \mathcal{N}\backslash\{1\} \label{cs1_nodal_blending} \\
& \;\; \sum\limits_{i \in \mathcal{I}}F_{ijn} = d_{jn} \quad \forall \, j \in \mathcal{J}, n \in \mathcal{N}\backslash\{1\} \label{cs1_nodal_demand} \\
& \;\; \sum\limits_{j \in \mathcal{J}}F_{ijn} \leq C_{i} \sum\limits_{n' \in \mathcal{P}_{n}}x_{in'} \quad \forall \, i \in \mathcal{I}, n \in \mathcal{N}\backslash\{1\} \label{cs1_nodal_capacity1} \\
& \;\; x_{in} + \sum\limits_{n' \in \mathcal{P}_{n}}x_{in'} \leq 1 \quad \forall \, i \in \mathcal{I}, n \in \mathcal{N}\backslash\{1\cup\mathcal{N}^{\mathrm{L}} \} \label{cs1_nodal_capacity2}\\
& \;\; 0 \leq c_{jn} \leq 1 \quad \forall \, j \in \mathcal{J}, n \in \mathcal{N}\backslash\{1\} \label{cs1_nodal_bound_c}\\
& \;\; d_{jn}^{\mathrm{min}} \leq d_{jn} \leq d_{jn}^{\mathrm{max}} \quad \forall \, j \in \mathcal{J}, n \in \mathcal{N}\backslash\{1\} \label{cs1_nodal_bound_d}\\
& \;\; F_{ijn} \geq 0 \quad \forall \, i \in \mathcal{I}, j \in \mathcal{J}, n \in \mathcal{N}\backslash\{1\} \label{cs1_nodal_bound_F}\\
& \;\; x_{in} \in \{0,1\} \quad \forall \, i \in \mathcal{I}, n \in \mathcal{N}\backslash\mathcal{N}^{\mathrm{L}}, \label{cs1_nodal_bound_x} 
\end{align}
where $\mathcal{N}$ and $\mathcal{N}^{\mathrm{L}}$ are the sets representing all nodes and leaf nodes in the scenario tree, respectively. Set $\mathcal{P}_{n}$ denotes the path from the root node to the parent node of node $n$ in the scenario tree. The stage containing node $n$ is denoted by $t(n)$. The probability of occurrence of node $n$ is denoted by $p_{n}$.

To be able to apply column generation, we introduce a set of auxiliary binary variables, $w_{in}$, defined by constraints \eqref{cs1_modified_1}, which help transform the above formulation to the same form as (MSSP). Further, constraints \eqref{cs1_nodal_capacity1} and \eqref{cs1_nodal_capacity2} are replaced by constraints \eqref{cs1_modified_2} and \eqref{cs1_modified_3}, respectively. 
\begin{align}
& \;\; w_{in} = \sum\limits_{n' \in \mathcal{P}_n}x_{in'} \quad \forall \, i \in \mathcal{I}, n \in \mathcal{N}\backslash\{1\} \label{cs1_modified_1} \\
& \;\; \sum\limits_{j \in \mathcal{J}}F_{ijn} \leq C_{i} w_{in} \quad \forall \, i \in \mathcal{I}, n \in \mathcal{N}\backslash\{1\} \label{cs1_modified_2} \\
& \;\; w_{in} + x_{in} \leq 1 \quad \forall \, i \in \mathcal{I}, n \in \mathcal{N}\backslash\{1 \cup\mathcal{N}^{\mathrm{L}} \} \label{cs1_modified_3} \\
& \;\; w_{in} \in \{0,1\} \quad \forall \, i \in \mathcal{I}, n \in \mathcal{N}\backslash\{1\}. \label{cs1_modified_4}
\end{align}
Note that although variable $w_{in}$ corresponds to node $n$ in the scenario tree, from constraints \eqref{cs1_modified_1}, it is easy to see that physically it points at whether an input tank $i$ has been installed in any of the previous stages (or equivalently up to the parent of node $n$). The major benefit of this is that when these variables are priced out as columns, sharing them among sibling nodes will ensure nonanticipativity (because they point to decisions in previous stages).

\subsection*{B. Column generation for the multistage blending case study}

The master problem for the column generation algorithm is as follows:
\begin{align}
    \maximize\limits_{\rho,x} {}&\quad \sum_{n \in \mathcal{N}\backslash\{1\}}\sum_{k \in \mathcal{K}_{n}}\rho_{nk}\Phi^{*}_{nk}-\sum_{n \in \mathcal{N}\backslash\mathcal{N}^{\mathrm{L}}}p_{n}\sum_{i \in \mathcal{I}}q_{i,t(n)}x_{in}\\
    \sto & \;\; \sum_{k \in \mathcal{K}_{n}}\rho_{nk}w_{ink}^{*}+x_{in} \leq 1 \quad \forall \, i \in \mathcal{I}, n \in \mathcal{N}\backslash\{1\cup\mathcal{N}^{\mathrm{L}}\} \qquad \quad \,\, [\gamma_{in}^{(1)}] \label{cs1_mp_cons1}\\
    & \;\; \sum_{k \in \mathcal{K}_{n}}\rho_{nk}w_{ink}^{*} = \sum_{n' \in \mathcal{P}_{n}}x_{in'} \quad \forall \, i \in \mathcal{I}, n \in \mathcal{N}\backslash\{1\} \qquad \qquad \quad [\gamma_{in}^{(2)}] \label{cs1_mp_cons2} \\
    & \;\; \sum_{k \in \mathcal{K}_{n}} \rho_{nk} = 1 \quad \forall \, n \in \mathcal{N}\backslash\{1\} \qquad \qquad \qquad  \qquad \qquad \qquad \quad \, [\mu_{n}] \label{cs1_convexity}\\
    & \;\; \rho_{nk} \in \{0,1\} \quad \forall \, n \in \mathcal{N}\backslash\{1\}, k \in \mathcal{K}_{n} \\
    & \;\; x_{in} \in \{0,1\} \quad \forall \, i \in \mathcal{I}, n \in \mathcal{N}\backslash\mathcal{N}^{\mathrm{L}},
\end{align}
where $\mathcal{K}_{n}$ is the set of columns priced out by the subproblem for node $n$ and binary variable $\rho_{nk}$ equals 1 if column $k$ is selected from $\mathcal{K}_{n}$. Convexity constraints \eqref{cs1_convexity} ensure that exactly one column is selected from the set $\mathcal{K}_{n}$. The dual prices corresponding to constraints \eqref{cs1_mp_cons1}, \eqref{cs1_mp_cons2}, and \eqref{cs1_convexity} are denoted by $\gamma^{(1)}$, $\gamma^{(2)}$, and $\mu$, respectively. Lastly, $\Phi_{nk}^{*}$ denotes the cost of column $k$ from pricing problem $n$ and is defined as follows: 
\begin{align}
\Phi_{nk|n>1}^{*} = \max_{c,d,F} {}&\quad p_{n}\Big[\sum_{i \in \mathcal{I}}\sum_{j \in \mathcal{J}}f_{j,t(n)}(c_{jn}) F_{ijn} - \sum_{i \in \mathcal{I}}\sum_{j \in \mathcal{J}}r_{ij}F_{ijn} - \sum_{i \in \mathcal{I}}b_{i,t(n)}\sum_{j \in \mathcal{J}}F_{ijn}\Big] &\\
\sto & \;\; \sum\limits_{i \in \mathcal{I}}\lambda_{i}F_{ijn} = c_{jn}d_{jn} \quad \forall \, j \in \mathcal{J} \label{cs1_colcost_1}\\
& \;\; \sum\limits_{i \in \mathcal{I}}F_{ijn} = d_{jn} \quad \forall \, j \in \mathcal{J} \label{cs1_colcost_2}\\
& \;\; \sum\limits_{j \in \mathcal{J}}F_{ijn} \leq C_{i} w_{ink}^{*} \quad \forall \, i \in \mathcal{I} \label{cs1_colcost_3}\\
& \;\; 0 \leq c_{jn} \leq 1 \quad \forall \, j \in \mathcal{J} \label{cs1_colcost_4}\\
& \;\; d_{jn}^{\mathrm{min}} \leq d_{jn} \leq d_{jn}^{\mathrm{max}} \quad \forall \, j \in \mathcal{J} \label{cs1_colcost_5}\\
& \;\; F_{ijn} \geq 0 \quad \forall \, i \in \mathcal{I}, j \in \mathcal{J}. \label{cs1_colcost_6}
\end{align}

Let $\Omega_{n} := \{F_{n},c_{n},d_{n}:\eqref{cs1_colcost_1}, \eqref{cs1_colcost_2}, \eqref{cs1_colcost_4}-\eqref{cs1_colcost_6}\}$. The pricing problem for node $n$ is defined as follows:
\begin{align}
\maximize_{w,c,d,F} {}&\quad p_{n}\Big[\sum_{i \in \mathcal{I}}\sum_{j \in \mathcal{J}}f_{j,t(n)}(c_{jn}) F_{ijn} - \sum_{i \in \mathcal{I}}\sum_{j \in \mathcal{J}}r_{ij}F_{ijn} \\ 
& \qquad \qquad - \sum_{i \in \mathcal{I}}b_{i,t(n)}\sum_{j \in \mathcal{J}}F_{ijn}\Big] - \sum_{i \in \mathcal{I}} (\gamma_{in}^{(1)} + \gamma_{in}^{(2)})w_{in} -  \mu_{n}  \nonumber \\
\sto & \;\; \sum\limits_{j \in \mathcal{J}}F_{ijn} \leq C_{i} w_{in} \quad \forall \, i \in \mathcal{I} \\
& \;\; w_{in} \in \{0,1\} \quad \forall \, i \in \mathcal{I} \\
& \;\; (F_{n}, c_{n}, d_{n}) \in \Omega_{n}.
\end{align}
Note that there is no pricing problem associated with the root node because $w_{in}$ is not defined at the root node.

\subsection*{C. Stochastic formulation for the mobile generator routing case study}

The stochasticity in the load at each bus of the power distribution network is modeled via the load variation factor, i.e $\beta_{it}$. The decision sequence in a time period is as follows: the decision to move gensets to new locations is made at the beginning of the time period, followed by realization of uncertainty in power demand at each bus, and finally, all operational and scheduling decisions are made. The deterministic formulation in Section \ref{sec:OPF} is extended to the stochastic programming formulation as follows: 
\begin{align}
\minimize {}&\quad \sum\limits_{n \in \mathcal{N}\backslash\{1 \cup \mathcal{N}^\mathrm{L}\}} p_{n} \sum_{m \in \mathcal{M}} C^{\mathrm{t}}_{m,t(n)}z_{mn} + \\
& \quad \sum\limits_{n \in \mathcal{N}\backslash\{1\}} p_{n}\alpha_{t(n)} \sum_{i \in \mathcal{B}} \sum_{h \in \mathcal{H}}\bigg[ C_{i,t(n),h}^{\mathrm{p}}p_{inh}^{\mathrm{g}}+  C^{\mathrm{VoLL}}_i (1-\delta_{inh}) p^{\mathrm{d}}_{inh} + \sum_{m \in \mathcal{M}}{C^{\mathrm{a}}_{t(n)}p_{minh}^{\mathrm{a}}}\bigg] \nonumber \\
\sto & \quad \sum_{i \in \mathcal{B}}y_{min}=1 \quad \forall \, m \in \mathcal{M}, n \in \mathcal{N}\backslash \mathcal{N}^\mathrm{L} \label{cs2_stoch_location} \\
& \quad y_{min}-y_{m,i,a(n)} \leq z_{mn} \quad \forall \, m \in \mathcal{M}, i \in \mathcal{B}, n \in \mathcal{N}\backslash\{1 \cup \mathcal{N}^\mathrm{L}\} \label{cs2_stoch_log1}\\
& \quad z_{mn} \leq 2 - y_{min} - y_{m,i,a(n)} \quad \forall \, m \in \mathcal{M}, i \in \mathcal{B}, n \in \mathcal{N}\backslash\{1 \cup \mathcal{N}^\mathrm{L}\} \label{cs2_stoch_log2}\\
& \quad p_{minh}^{\mathrm{a}} \leq \overline{p}_{m}^{\mathrm{a}} y_{m,i,a(n)} \quad \forall \, m \in \mathcal{M}, i \in \mathcal{B}, n \in \mathcal{N}\backslash\{1\}, h \in \mathcal{H} \label{cs2_stoch_active_power_bound} \\
& \quad  -p_{minh}^{\mathrm{a}} \leq q_{minh}^{\mathrm{a}} \leq p_{minh}^{\mathrm{a}} \quad \forall \, m \in \mathcal{M}, i \in \mathcal{B}, n \in \mathcal{N}\backslash\{1\}, h \in \mathcal{H} \label{cs2_stoch_reactive_power_bounds} \\
& \quad p_{inh}^{\mathrm{g}} - {p_{inh}^{\mathrm{d}}} \delta_{inh} + \sum_{m \in \mathcal{M}}p_{minh}^{\mathrm{a}} = \nonumber \\ 
& \quad \qquad \qquad \sum_{j \in \mathcal{C}(i)} P_{jnh} - \sum_{j \in \mathcal{D}(i)}{(P_{inh}-r_{i}l_{inh})} + g_{i}v_{inh} \quad \forall \, i \in \mathcal{B}, n \in \mathcal{N}\backslash\{1\}, h \in \mathcal{H} \label{cs2_stoch_powerflow_1} \\
& \quad q_{ith}^{\mathrm{g}} - {q_{inh}^{\mathrm{d}}} \delta_{inh} + \sum_{k \in \mathcal{K}}q_{minh}^{\mathrm{a}} = \nonumber \\
& \quad \qquad \qquad \sum_{j \in \mathcal{C}(i)} Q_{jnh} - \sum_{j \in \mathcal{D}(i)}{(Q_{inh}-x_{i}l_{inh})} + b_{i}v_{inh} \quad \forall \, i \in \mathcal{B}, n \in \mathcal{N}\backslash\{1\}, h \in \mathcal{H} \label{cs2_stoch_powerflow_2} \\
& \quad v_{inh} = v_{d(i),n,h} - 2(r_{i}P_{inh} + x_{i}Q_{inh}) + (r_{i}^{2} + x_{i}^{2})l_{inh} \quad \forall \, i \in \mathcal{L}, n \in \mathcal{N}\backslash\{1\}, h \in \mathcal{H} \label{cs2_stoch_powerflow_3}\\
& \quad l_{inh}v_{d(i),n,h} \geq P_{inh}^2 + Q_{inh}^2 \quad \forall \, i \in \mathcal{L}, n \in \mathcal{N}\backslash\{1\}, h \in \mathcal{H} \label{cs2_stoch_powerflow_4}\\
& \quad P_{inh}^2 + Q_{inh}^2 \leq A_{i}^2 \quad \forall \, i \in \mathcal{L}, n \in \mathcal{N}\backslash\{1\}, h \in \mathcal{H} \label{cs2_stoch_powerflow_5}\\
& \quad (P_{inh}-r_{i}l_{inh})^2 + (Q_{inh}-x_{i}l_{inh})^2 \leq A_{i}^2 \quad \forall \, i \in \mathcal{L}, n \in \mathcal{N}\backslash\{1\}, h \in \mathcal{H} \label{cs2__stoch_powerflow_6}\\
& \quad 0 \leq p_{inh}^{\mathrm{g}} \leq \overline{p}_{i}^{\mathrm{g}} \quad \forall \, i \in \mathcal{B}, n \in \mathcal{N}\backslash\{1\}, h \in \mathcal{H} \label{cs2_stoch_bound_pg}\\
& \quad 0 \leq q_{inh}^{\mathrm{g}} \leq \overline{q}_{i}^{\mathrm{g}} \quad \forall \, i \in \mathcal{B}, n \in \mathcal{N}\backslash\{1\}, h \in \mathcal{H} \label{cs2_stoch_bound_qg}\\
& \quad  0 \leq \delta_{inh} \leq 1 \quad \forall \, i \in \mathcal{B}, n \in \mathcal{N}\backslash\{1\}, h \in \mathcal{H} \label{cs2_stoch_bound_delta}\\
& \quad \underline{v}_{i} \leq v_{inh} \leq \overline{v}_{i} \quad \forall \, i \in \mathcal{B}, n \in \mathcal{N}\backslash\{1\}, h \in \mathcal{H} \label{cs2_stoch_bound_v}\\
& \quad l_{inh} \geq 0 \quad \forall \, i \in \mathcal{L}, n \in \mathcal{N}\backslash\{1\}, h \in \mathcal{H} \label{cs2_stoch_bound_l}\\
& \quad  p_{minh}^{\mathrm{a}} \geq 0 \quad \forall \, m \in \mathcal{M}, i \in \mathcal{B}, n \in \mathcal{N}\backslash\{1\}, h \in \mathcal{H} \label{cs2_stoch_bound_pa}\\
& \quad y_{min} \in \{0,1\} \quad \forall \, m \in \mathcal{M}, i \in \mathcal{B}, n \in \mathcal{N}\backslash \mathcal{N}^\mathrm{L}  \label{cs2_stoch_bound_y}\\
& \quad z_{mn} \in \{0,1\} \quad \forall \, m \in \mathcal{M}, n \in \mathcal{N}\backslash\{1 \cup \mathcal{N}^\mathrm{L}\} \label{cs2_stoch_bound_z}\\
\operatorname{where} & \quad p_{inh}^{\mathrm{d}} = p^{\mathrm{d,base}}_{i} \beta_{in} e_{t(n),h} \quad \forall \, i \in \mathcal{B}, n \in \mathcal{N}\backslash\{1\}, h \in \mathcal{H} \\
& \quad q_{inh}^{\mathrm{d}} = q^{\mathrm{d,base}}_i \beta_{in} e_{t(n),h} \quad \forall \, i \in \mathcal{B}, n \in \mathcal{N}\backslash\{1\}, h \in \mathcal{H}.
\end{align}
The sets $\mathcal{N}$ and $\mathcal{N}^{\mathrm{L}}$ represent all and leaf nodes in the scenario tree, respectively. The parent node of a node $n$ is denoted by $a(n)$. The stage containing node $n$ is denoted by $t(n)$. The probability of occurrence of node $n$ is denoted by $p_{n}$.

Although the above formulation fits the structure defined by (MSSP), for column sharing to produce meaningful columns that satisfy the nonanticipativity constraints, we introduce an auxiliary binary variable $w_{min}$, defined by constraints \eqref{cs2_modified_1}. Although $w_{min}$ will be part of the subproblem corresponding to node $n$, it actually represents the decision made on the parent node of node $n$. Therefore, sharing it with siblings of node $n$ will ensure nonanticipativity. Also, constraints \eqref{cs2_stoch_location}-\eqref{cs2_stoch_active_power_bound} are replaced by \eqref{cs2_modified_2}-\eqref{cs2_modified_5}. 
\begin{align}
& \;\; w_{min} = y_{m,i,a(n)} \quad \forall \, m \in \mathcal{M}, i \in \mathcal{B}, n \in \mathcal{N}\backslash\{1\} \label{cs2_modified_1} \\
& \;\; \sum_{i \in \mathcal{B}}w_{min} = 1 \quad \forall \, m \in \mathcal{M}, n \in \mathcal{N}\backslash\{1\} \label{cs2_modified_2} \\
& \;\; y_{min}-w_{min} \leq z_{mn} \quad \forall \, m \in \mathcal{M}, i \in \mathcal{B}, n \in \mathcal{N}\backslash\{1 \cup \mathcal{N}^\mathrm{L}\} \label{cs2_modified_3}\\
& \;\; z_{mn} \leq 2 - y_{min} - w_{min} \quad \forall \, m \in \mathcal{M}, i \in \mathcal{B}, n \in \mathcal{N}\backslash\{1 \cup \mathcal{N}^\mathrm{L}\} \label{cs2_modified_4}\\
& \;\; p_{minh}^{\mathrm{a}} \leq \overline{p}_{m}^{\mathrm{a}} w_{min} \quad \forall \, m \in \mathcal{M}, i \in \mathcal{B}, n \in \mathcal{N}\backslash\{1\}, h \in \mathcal{H} \label{cs2_modified_5} \\
& \;\; w_{min} \in \{0,1\} \quad \forall \, m \in \mathcal{M}, i \in \mathcal{B}, n \in \mathcal{N}\backslash\{1\}. \label{cs2_modified_6}
\end{align}

\subsection*{D. Column generation decomposition for the mobile generator routing case study}

The master problem for the column generation algorithm is as follows:
\begin{align}
    \minimize\limits_{\rho,y,z} {}&\quad \sum_{n \in \mathcal{N}\backslash\{1\}}\sum_{\mathcal{K}_{n}}\rho_{nk}\Phi^{*}_{nk}+\sum_{n \in \mathcal{N}\backslash\{1 \cup \mathcal{N}^{\mathrm{L}}\}}p_{n}\sum_{m \in \mathcal{M}}C^{\mathrm{t}}_{m,t(n)}z_{mn}\\
    \sto  & \;\; \sum_{k \in \mathcal{K}_{n}}\rho_{nk}w_{mink}^{*} = y_{m,i,a(n)} \quad \forall \, m \in \mathcal{M}, i \in \mathcal{B}, n \in \mathcal{N}\backslash\{1\} \qquad \qquad \qquad \;\;\; [\gamma^{(1)}_{min}] \label{cs2_mp_cons1} \\
    & \;\; \sum_{k \in \mathcal{K}_{n}}\rho_{nk}w_{mink}^{*} - y_{min} + z_{mn} \geq 0 \quad \forall \, m \in \mathcal{M}, i \in \mathcal{B}, n \in \mathcal{N}\backslash\{1 \cup \mathcal{N}^{\mathrm{L}}\} \quad [\gamma^{(2)}_{min}] \label{cs2_mp_cons2}\\
    & \;\; z_{mn} + y_{min} + \sum_{k \in \mathcal{K}_{n}}\rho_{nk}w_{mink}^{*} \leq 2 \quad \forall \, m \in \mathcal{M}, i \in \mathcal{B}, n \in \mathcal{N}\backslash\{1 \cup \mathcal{N}^{\mathrm{L}}\} \, \quad [\gamma^{(3)}_{min}] \label{cs2_mp_cons3}\\
    & \;\; \sum_{k \in \mathcal{K}_{n}} \rho_{nk} = 1 \quad \forall \, n \in \mathcal{N}\backslash\{1\} \qquad \qquad \qquad  \qquad \qquad \qquad \qquad \qquad \qquad \quad \;\; [\mu_{n}] \label{cs2_convexity}\\
    & \;\; \rho_{nk} \in \{0,1\} \quad \forall \, n \in \mathcal{N}\backslash\{1\}, k \in \mathcal{K}_{n} \\
    & \;\; y_{min} \in \{0,1\} \quad \forall \, m \in \mathcal{M}, i \in \mathcal{B}, n \in \mathcal{N}\backslash\mathcal{N}^{\mathrm{L}} \\
    & \;\; z_{mn} \in \{0,1\} \quad \forall \, m \in \mathcal{M}, n \in \mathcal{N}\backslash\{1 \cup \mathcal{N}^{\mathrm{L}}\},
\end{align}
where $\mathcal{K}_{n}$ is the set of columns priced out by the subproblem for node $n$ and binary variable $\rho_{nk}$ equals 1 if column $k$ is selected from $\mathcal{K}_{n}$. Convexity constraints \eqref{cs2_convexity} ensure that exactly one column gets selected from the set $\mathcal{K}_{n}$. The dual prices corresponding to constraints \eqref{cs2_mp_cons1}, \eqref{cs2_mp_cons2}, \eqref{cs2_mp_cons3}, and \eqref{cs2_convexity} are denoted by $\gamma^{(1)}$, $\gamma^{(2)}$, $\gamma^{(3)}$, and $\mu$, respectively. Lastly, $\Phi_{nk}^{*}$ denotes the cost of column $k$ from pricing problem $n$ and is defined as follows: 
\begin{align}
\Phi_{nk|n>1}^{*} = \min {}&\quad p_{n}\alpha_{t(n)} \sum_{i \in \mathcal{B}} \sum_{h \in \mathcal{H}}\bigg[ C_{i,t(n),h}^{\mathrm{p}}p_{inh}^{\mathrm{g}}+  C^{\mathrm{VoLL}}_i (1-\delta_{inh}) p^{\mathrm{d}}_{inh} + \sum_{m \in \mathcal{M}}{C^{\mathrm{a}}_{t(n)}p_{minh}^{\mathrm{a}}}\bigg] &\\
\sto & \quad p_{minh}^{\mathrm{a}} \leq \overline{p}_{m}^{\mathrm{a}} w_{mink}^{*} \quad \forall \, m \in \mathcal{M}, i \in \mathcal{B}, h \in \mathcal{H} \\
& \quad  -p_{minh}^{\mathrm{a}} \leq q_{minh}^{\mathrm{a}} \leq p_{minh}^{\mathrm{a}} \quad \forall \, m \in \mathcal{M}, i \in \mathcal{B}, h \in \mathcal{H} \label{cs2_colcost_reactivepowergen_bound} \\
& \quad p_{inh}^{\mathrm{g}} - {p_{inh}^{\mathrm{d}}} \delta_{inh} + \sum_{m \in \mathcal{M}}p_{minh}^{\mathrm{a}} = \nonumber \\ 
& \quad \qquad \qquad \sum_{j \in \mathcal{C}(i)} P_{jnh} - \sum_{j \in \mathcal{D}(i)}{(P_{inh}-r_{i}l_{inh})} + g_{i}v_{inh} \quad \forall \, i \in \mathcal{B}, h \in \mathcal{H} \\
& \quad q_{ith}^{\mathrm{g}} - {q_{inh}^{\mathrm{d}}} \delta_{inh} + \sum_{k \in \mathcal{K}}q_{minh}^{\mathrm{a}} = \nonumber \\
& \quad \qquad \qquad \sum_{j \in \mathcal{C}(i)} Q_{jnh} - \sum_{j \in \mathcal{D}(i)}{(Q_{inh}-x_{i}l_{inh})} + b_{i}v_{inh} \quad \forall \, i \in \mathcal{B}, h \in \mathcal{H} \\
& \quad v_{inh} = v_{d(i),n,h} - 2(r_{i}P_{inh} + x_{i}Q_{inh}) + (r_{i}^{2} + x_{i}^{2})l_{inh} \quad \forall \, i \in \mathcal{L}, h \in \mathcal{H} \\
& \quad l_{inh}v_{d(i),n,h} \geq P_{inh}^2 + Q_{inh}^2 \quad \forall \, i \in \mathcal{L}, h \in \mathcal{H} \\
& \quad P_{inh}^2 + Q_{inh}^2 \leq A_{i}^2 \quad \forall \, i \in \mathcal{L}, h \in \mathcal{H} \\
& \quad (P_{inh}-r_{i}l_{inh})^2 + (Q_{inh}-x_{i}l_{inh})^2 \leq A_{i}^2 \quad \forall \, i \in \mathcal{L}, h \in \mathcal{H} \\
& \quad 0 \leq p_{inh}^{\mathrm{g}} \leq \overline{p}_{i}^{\mathrm{g}} \quad \forall \, i \in \mathcal{B}, h \in \mathcal{H} \\
& \quad 0 \leq q_{inh}^{\mathrm{g}} \leq \overline{q}_{i}^{\mathrm{g}} \quad \forall \, i \in \mathcal{B}, h \in \mathcal{H} \\
& \quad  0 \leq \delta_{inh} \leq 1 \quad \forall \, i \in \mathcal{B}, h \in \mathcal{H} \\
& \quad \underline{v}_{i} \leq v_{inh} \leq \overline{v}_{i} \quad \forall \, i \in \mathcal{B}, h \in \mathcal{H} \\
& \quad l_{inh} \geq 0 \quad \forall \, i \in \mathcal{L}, h \in \mathcal{H} \\
& \quad  p_{minh}^{\mathrm{a}} \geq 0 \quad \forall \, m \in \mathcal{M}, i \in \mathcal{B}, h \in \mathcal{H}. \label{cs2_colcost_bound_pa}
\end{align}

Let $\Omega_{n} := \{p^{g}_{n},q^{g}_{n},p^{a}_{n},q^{a}_{n},P_{n},Q_{n},\delta_{n},v_{n},l_{n}:\eqref{cs2_colcost_reactivepowergen_bound}-\eqref{cs2_colcost_bound_pa}\}$, then the pricing problem for node $n$ is defined as follows:
\begin{align} 
\minimize {}&\quad p_{n}\alpha_{t(n)} \sum_{i \in \mathcal{B}} \sum_{h \in \mathcal{H}}\bigg[ C_{i,t(n),h}^{\mathrm{p}}p_{inh}^{\mathrm{g}}+  C^{\mathrm{VoLL}}_i (1-\delta_{inh}) p^{\mathrm{d}}_{inh} + \sum_{m \in \mathcal{M}}{C^{\mathrm{a}}_{t(n)}p_{minh}^{\mathrm{a}}}\bigg] - \\
& \qquad \qquad \qquad \qquad \qquad \qquad \qquad  \sum_{m \in \mathcal{M}}\sum_{i \in \mathcal{B}} (\gamma^{(1)}_{min}+\gamma^{(2)}_{min}+\gamma^{(3)}_{min})w_{min} -  \mu_{n}  \nonumber \\
\sto & \quad \sum_{i \in \mathcal{B}} w_{min} = 1 \quad \forall \, m \in \mathcal{M} \\ 
& \quad  p_{minh}^{\mathrm{a}} \leq \overline{p}_{m}^{\mathrm{a}} w_{min} \quad \forall \, m \in \mathcal{M}, i \in \mathcal{B}, h \in \mathcal{H}  \\
& \quad w_{min} \in \{0,1\} \quad \forall \, m \in \mathcal{M}, i \in \mathcal{B} \\ 
& \quad (p^{g}_{n},q^{g}_{n},p^{a}_{n},q^{a}_{n},P_{n},Q_{n},\delta_{n},v_{n},l_{n}) \in \Omega_{n}.
\end{align}
Note that there is no subproblem for the root node because $w_{min}$ is not defined at the root node.

\subsection*{E. Parameter distributions for case studies}

\begin{table}[htbp]
  \centering
    \begin{tabular}{lc}
    \toprule
    \textbf{Parameter} & 
    \textbf{Distribution}\\
    \midrule
    Input tank capacity, $C_{i}$ & $\mathcal{U}(0.5,3)10^3$   \\
    Input tank quality spec., $\lambda_{i}$ & $\mathcal{U}(0,1)$\\
    Installation cost, $q_{it}$ & $(|\mathcal{T}|-t+1)\mathcal{U}(200,400)10^3C_{i}/\sum_{i \in \mathcal{I}}C_{i}$ \\
    Production cost, $b_{it}$ & $\mathcal{U}(10,20)$ \\
    Transport cost, $r_{ij}$  & $\mathcal{U}(0.1,0.2)$   \\
    Minimum demand, $d_{jt}^{\mathrm{min}}$ & $\mathcal{U}(0.1,2)10^3$ \\
    Maximum demand, $d_{jt}^{\mathrm{max}}$ &  $\mathcal{U}(5,10)10^3$ \\
    Price slope, $m$         & $\mathcal{U}(20,50)$\\
    Price intercept, $l$     & $\mathcal{U}(5,20)$\\
    \bottomrule
    \end{tabular}%
\caption{Distributions used for generating parameters for instances in the multistage blending case study.}
  \label{tab:ParDistributionTable_blending}%
\end{table}%

\begin{table}[htbp]
  \centering
    \begin{tabular}{ll}
    \toprule
    \textbf{Parameter} & 
    \textbf{Distribution}\\
    \midrule
    Power gen. capacity of genset $k$, $\overline{p}_{k}^{\mathrm{a}}$ & $\mathcal{U}(0.1,0.8)$ \\
   Variable generation cost of gensets in time period $t$, $C^{\mathrm{a}}_{t}$ & $\mathcal{U}(100,120)$ \\
    Cost of moving genset $m$ in time period $t$, $C^{\mathrm{t}}_{mt}$ & $\mathcal{U}(25,35)$ \\
    Variable generation cost of stationary units at bus $i$ at time point $h$, $C_{ith}^{\mathrm{p}}$ & $C_{i}^{\mathrm{p,base}} \ \mathcal{U}(0.8,1.2)$\\
    Load variation factor for bus $i$ in time period $t$, $\beta_{it}$ & $\mathcal{U}(0.5,1.5)$\\
    Time of day factor at hour $h$, $e_{th}$ & $\mathcal{U}(0.5,1)$\\
    Upper bound on active power gen. at bus $i$, $\overline{p}_{i}^{\mathrm{g}}$ &  $\overline{p}_{i}^{\mathrm{g,base}} \ \mathcal{U}(0.8,1.2)$ \\
    Upper bound on reactive power gen. at bus $i$, $\overline{q}_{i}^{\mathrm{g}}$ &  $\overline{q}_{i}^{\mathrm{g,base}} \ \mathcal{U}(0.8,1.2)$ \\

    \bottomrule
    \end{tabular}%
\caption{Distributions used for generating parameters for instances in the power distribution network case study. $\overline{p}_{k}^{\mathrm{a}}$,  $\overline{p}_{i}^{\mathrm{g}}$, and $\overline{q}_{i}^{\mathrm{g}}$ are in units of MW. $C^{\mathrm{a}}_{t}$ and $C_{ith}^{\mathrm{p}}$ are in units of \$/MWh. $C^{\mathrm{t}}_{mt}$ is in units of \$. }
  \label{tab:ParDistributionTable_opf}%
\end{table}%

\begin{table}[htbp]
\fontsize{9}{16}\selectfont
\setlength\tabcolsep{5pt}
  \centering
\begin{tabular}{crrrrrrrrrrrrr}
\hline
Bus/line $i$  & $r_i$      & $x_i$      & $A_i$     & $p_{i}^{\mathrm{d,base}}$ & $q_{i}^{\mathrm{d,base}}$ & $b_i$       & $g_i$ & $C_{i}^{\mathrm{p,base}}$ & $C_{i}^
{\mathrm{VoLL}}$ & $\overline{p}_{i}^{\mathrm{g,base}}$  & $\overline{q}_{i}^{\mathrm{g,base}}$  & $\underline{v}_{i}$ & $\overline{v}_{i}$ \\ \hline
0  &        &        &       & 0       & 0       & 0       & 0 & 50      & 3500  & 1$\times$10\textsuperscript{6} & 1$\times$10\textsuperscript{6} & 1        & 1        \\
1  & 0.001  & 0.12   & 2     & 0.3872  & 0.0910  & -0.0011 & 0 & 0       & 3500  & 0        & 0        & 0.81     & 1.21     \\
2  & 0.0883 & 0.1262 & 0.256 & 0       & 0       & -0.0028 & 0 & 0       & 3500  & 0        & 0        & 0.81     & 1.21     \\
3  & 0.1384 & 0.1978 & 0.256 & 0.2402  & 0.0568  & -0.0024 & 0 & 0       & 3500  & 0        & 0        & 0.81     & 1.21     \\
4  & 0.0191 & 0.0273 & 0.256 & 0.2346  & 0.0486  & -0.0004 & 0 & 0       & 3500  & 0        & 0        & 0.81     & 1.21     \\
5  & 0.0175 & 0.0251 & 0.256 & 0.2582  & 0.0546  & -0.0008 & 0 & 0       & 3500  & 0        & 0        & 0.81     & 1.21     \\
6  & 0.0482 & 0.0689 & 0.256 & 0.2438  & 0.0510  & -0.0006 & 0 & 0       & 3500  & 0        & 0        & 0.81     & 1.21     \\
7  & 0.0523 & 0.0747 & 0.256 & 0       & 0.0438  & -0.0006 & 0 & 0       & 3500  & 0.1969   & 0        & 0.81     & 1.21     \\
8  & 0.0407 & 0.0582 & 0.256 & 0.247   & 0.0518  & -0.0012 & 0 & 0       & 3500  & 0        & 0        & 0.81     & 1.21     \\
9  & 0.01   & 0.0143 & 0.256 & 0.2458  & 0.0684  & -0.0004 & 0 & 0       & 3500  & 0        & 0        & 0.81     & 1.21     \\
10 & 0.0241 & 0.0345 & 0.256 & 0.2434  & 0.0530  & -0.0004 & 0 & 0       & 3500  & 0        & 0        & 0.81     & 1.21     \\
11 & 0.0103 & 0.0148 & 0.256 & 0.2264  & 0.0466  & -0.0001 & 0 & 35 & 3500  & 0.4   & 0.4   & 0.81     & 1.21     \\
12 & 0.001  & 0.12   & 1     & 0.4438  & 0.0982  & -0.0001 & 0 & 0       & 3500  & 0        & 0        & 0.81     & 1.21     \\
13 & 0.1559 & 0.1119 & 0.204 & 0.2028  & 0.0416  & -0.0002 & 0 & 0       & 3500  & 0        & 0        & 0.81     & 1.21     \\
14 & 0.0953 & 0.0684 & 0.204 & 0.2448  & 0.0566  & -0.0001 & 0 & 0       & 3500  & 0        & 0        & 0.81     & 1.21     \\ \hline
\end{tabular}
\caption{Bus and line related data for the optimal power flow case study. Quantities are given with power in units of MW and energy in units of MWh. A per-unit system is used for the distribution network with a base power of 1 MW, meaning that parameters such as $r_i$ and $x_i$ are given in dimensionless per-unit quantities. $C_{i}^{\mathrm{p,base}}$ and $C_{i}^
{\mathrm{VoLL}}$ are in units of \$/MWh.}
  \label{tab:bus_line_data}%
\end{table}

\end{document}